\documentclass[11pt]{article}

\usepackage[margin=1.in]{geometry}
\usepackage{mathtools}
\usepackage{amssymb} 
\usepackage{enumitem}
\usepackage{amsthm}
\usepackage{amsopn}
\usepackage{bbm}
\usepackage[numbers]{natbib}
\usepackage{amsmath}
\usepackage{amsfonts}
\usepackage{amssymb}
\usepackage{graphicx} 
\usepackage{enumitem}
\DeclarePairedDelimiterX\Basics[1](){ #1}
\usepackage[utf8]{inputenc}
\usepackage{color}
\usepackage{caption}
\usepackage{subcaption}
\usepackage{algorithm,algorithmic}
\usepackage{enumitem}
\usepackage{pifont}

\usepackage{pifont}

\usepackage[T1]{fontenc}
\usepackage{etoolbox}

\usepackage{float}
\newtheorem{theorem}{Theorem}
\newtheorem{corollary}{Corollary}[theorem]

\newtheorem{remark}[theorem]{Remark}
\newtheorem{definition}[theorem]{Definition}

\numberwithin{equation}{section}

\newcommand{\innermid}{\nonscript\;\delimsize\vert\nonscript\;}
\newcommand{\activatebar}{%
  \begingroup\lccode`\~=`\|
  \lowercase{\endgroup\let~}\innermid 
  \mathcode`|=\string"8000
}

\makeatletter

\newlist{steps}{enumerate}{1}
\setlist[steps, 1]{label = Step \arabic*:}

\newcommand{\interior}[1]{%
  {\kern0pt#1}^{\mathrm{o}}%
}

\def\hatR{\widehat{R}}

\def\hatV{\widehat{\sigma}^2}
\def\hatM{\widehat{M}}
\def\hatK{\widehat{K}}

\def\R{\mathbb{R}}

\def\HH{\mathcal{H}}

\def\PP{\mathcal{P}}
\def\AA{\mathcal{A}}

\def\E{\mathbb{E}}
\def\P{\mathbb{P}}
\def\R{\mathbb{R}}

\begin{document}
\title{Entropy Regularization for Mean Field Games with Learning}
\author{
Xin Guo
\thanks{Department of Industrial Engineering and Operations Research, University of California, Berkeley, and Amazon Research, USA. Email: xinguo@berkeley.edu.}
\and Renyuan Xu
\thanks{Epstein Department of Industrial \& Systems Engineering, University of Southern California, US and Mathematical Institute, University of Oxford,  UK. Email: renyuanx@usc.edu.}
\and Thaleia Zariphopoulou
\thanks{Departments of Mathematics and IROM, The University of Texas at Austin, Austin,
USA,  and the Oxford-Man Institute, University of Oxford, Oxford, UK. Email:  zariphop@math.utexas.edu.}
}
\date{\today}
\maketitle


\begin{abstract}
Entropy regularization has been extensively adopted   to improve the efficiency, the stability, and the convergence of algorithms in reinforcement learning.
This paper analyzes  both quantitatively and qualitatively the impact of entropy regularization for Mean Field Games (MFGs) with learning in a finite time horizon.  Our study provides a theoretical justification that entropy regularization yields time-dependent policies and, furthermore, helps stabilizing and accelerating  convergence to the game equilibrium. 
In addition, this study leads to a policy-gradient algorithm with exploration in MFG. With this algorithm,  agents are able to  learn the optimal exploration scheduling, with stable and fast  
convergence to the game equilibrium.
\end{abstract}

\maketitle
\section{Introduction}
Reinforcement learning (RL) is one of the three fundamental machine learning paradigms,
alongside supervised learning and unsupervised learning.
RL is learning via trial and error, through interactions with an environment and possibly with other agents; in RL,  an agent takes an action and receives a reinforcement signal in terms of a numerical reward, which encodes the outcome of her action.  In order to maximize the accumulated reward over time, the agent  learns to select her actions  based on her past experiences (exploitation) and/or  by making new choices (exploration).

Exploration and exploitation are  the essence of RL.
Exploration provides opportunities to improve from  current sub-optimal solutions to the ultimate global optimal one, yet is time consuming and computationally expensive  as over-exploration   may impair the convergence to the optimal solution. Meanwhile, pure exploitation, i.e., myopically picking the current solution based solely on past experience, though easy to implement, tends to yield  sub-optimal global solutions.
 Therefore, an appropriate trade-off between exploration-exploitation  is crucial for  RL algorithm design to improve  the learning and the optimization procedure.

 \paragraph{Entropy regularization.}
 One common approach to balance  the exploration-exploitation in RL is to introduce entropy regularization  \cite{ALNS2019,LB2019,NJG2017}.
In  RL settings with more than one agent,  there are two major sources of  uncertainty: the unknown environment and the actions of the other agents.
Shannon entropy and cross-entropy are two natural choices for entropy regularization:
the former quantifies the information gain of exploring the environment while the latter measures the benefit from exploring the actions of other agents. This information-theoretic perspective of exploration has been well understood in single-agent RL; see for instance \cite{GLF2018, hazan2019,LB2019, NJG2017, RV2016}.

However, there is virtually no theoretical study on the role of entropy regularization in multi-agent RL (MARL), with the exception of \cite{AKS2020}.
    Indeed, most existing studies are  empirical, demonstrating  convergence improvement and variance reduction when entropy regularization is added.
 For instance,  \cite{HSSCL2018} showed via empirical analysis that  policy features can be learned directly from pure observations of  other agents and that the non-stationarity of the environment can  be reduced by adding cross-entropy;
 \cite{HS2002} applied the cross-entropy regularization to demonstrate the convergence of fictitious play in a discrete-time model with a finite number of agents while \cite{RDSF2018} used the cross-entropy loss to train the
prediction of  other agents' actions via observations of  their behavior.  
The only theoretical work so far can be found in \cite{AKS2020} in an infinite  horizon setting in which
 a regularized Q-learning algorithm for stationary discrete-time mean field games was proposed along with its convergence analysis.  
Still, the problem remains open for the  finite time horizon case, which arise often in many applications in operations research, supply-chain management and finance.

\paragraph{Optimal exploration scheduling.}
 Another  major challenge for both single-agent RL and MARL is exploration efficiency.
  In practice, there are various heuristic designs of explorations for MARL, including adding random noise in the parameter space \cite{plappert2017}, the approach of $\varepsilon$-greedy   policy \cite{wunder2010classes}, and the method with softmax \cite{iqbal2019actor}. However, there is no theoretical validation of these approaches.
 
Recently,
 time-invariant Gaussian exploration was  applied to single-agent RL (\cite{HCDSDA2016,MTR2019,WZZ2018}) and  time-dependent ``optimal exploration scheduling'' was  derived for single-agent mean-variance portfolio selection problem in \cite{WZ2019}.
 In these works,  the degree of exploration was characterized by the variance of the Gaussian distribution and the term ``optimal exploration scheduling'' was coined for the time-dependent variance  of the Gaussian distribution.
 
Exploration schemes are inherently
time-dependent, as it is  necessary to balance the  free exploration at the initial phase and the greedy control policy towards  terminal time. Yet, it seems that there is no existing work on analyzing such  time-dependent learning policies for MARL, neither empirical or theoretical.  

\paragraph{Model-free vs model-based approach for MFG with learning.}
There are two popular approaches {in single-agent reinforcement-learning} to handle unknown or partially known environments: the model-based approach and the model-free approach.  In the model-based paradigm, the agent is assumed to know  the model structure but has no access to the model parameters. In this case, the agent  estimates the unknown model parameters and, then, constructs a control policy based on the knowledge of the model \cite{aastrom2013adaptive,dean2019}. In the model-free paradigm, the agent learns the optimal policy {\it directly} via interacting with the system, without inferring the model parameters. Examples of model-free approach include policy gradient method \cite{bhandari2019} and actor critic method \cite{FYCW2019}. In practice, due to the lack of information on the actual system, model-based approach tends to suffer from model mis-specification \cite{bhandari2019}. On the other hand, the execution of the model-free algorithm does not rely on the assumptions of the model, thus is more robust against model mis-specification \cite{hambly2020policy}. 

Given the robustness of the model-free approach  and the additional complexity from the game interactions,  model-free approach appears more appropriate for MFG with learning  where the representative agent faces uncertainties  about both the unknown environment and the large population of strategic opponents. 


\paragraph{Our work.}
In this paper, we propose to study entropy regularization for MARL with a large population, namely, within the framework of  the mean field game (MFG). This transition from MARL to MFG with learning is critical to avoid the  curse of dimensionality in MARL. 

We analyze both quantitatively and qualitatively the impact of entropy regularization in MFG with learning in a finite time horizon.
 We adopt two different entropies: first, the Shannon entropy and, then, a combination of Shannon entropy and the cross-entropy, which we call the {\it enhanced entropy}.

 \begin{itemize}
 \item We derive explicit Nash equilibrium (NE) solutions  (Theorems \ref{thm:current_shannon} and \ref{thm:current_cross}) for a class of linear-quadratic (LQ) stochastic games.
 Our study provides a theoretical justification to the fact that entropy regularization yields time-dependent policies. Furthermore, it helps stabilizing and accelerating  convergence to the game equilibrium. 

\item  This theoretical study enables us to design a model-free policy-gradient algorithm for  MFG with learning.  
Under this algorithm,  agents are able to  learn efficiently the optimal exploration scheduling in an unknown environment and with  a large group of competing agents.
The convergence to the game equilibrium is stable and fast when appropriate exploration rates are chosen. 

\end{itemize}

\paragraph{Additional related works.}
Our algorithm is  inspired by the recent success of policy-gradient method for single-agent LQ regulators \cite{FGKM2018}. In addition, there is a concurrent work on the global convergence of policy gradient for MFG \cite{wang2020global}, yet without exploration. We   also mention recent works on two-agent zero-sum LQ games \cite{ZYB2019}, general-sum LQ games \cite{hambly2021policy}, and the LQ mean field control problem with common noise \cite{CLT2019}.

  \paragraph{Organization.} The rest of the paper is organized as follows. Section \ref{sec:setup} provides the mathematical framework for MFG with learning, Section \ref{sec:LQ} focuses on analyzing the impact of Shannon entropy and the enhanced entropy in a class of LQ games, and  Section \ref{sec:experiment} proposes a policy-gradient based algorithm with  entropy regularization,
 and provides its numerical  performance.

\section{Mathematical Formulation}\label{sec:setup}

We start with the mathematical formulation of the MFG with learning. 

\paragraph{Key ideas.}

There are several key components for the formulation.

The first component is the {\it aggregation idea} from the theory of MFG to address the  curse of dimensionality in MARL. Specifically, it is to consider $N$ agents, and assume that they are all  identical, indistinguishable and interchangeable, and that interactions among them are based on the {\it macroscopic information}, which is the empirical state distribution and action distribution of all agents.  This allows us to work instead with a representative agent $i$, her state $X^i_t$, her policy $\pi^i_t$ at time $t\in [0, T]$,  and her interaction with other agents through the macroscopic information. Since agent $i$  depends on other agents only through the empirical measure, we may then consider both the 
population state distribution and action distribution if such limits exist when  $N \rightarrow \infty$. Moreover, the subscript $i$ can be dropped and one can focus on a representative agent in this MFG formulation since all agents are assumed to be identical and indistinguishable.

The second component is  how to model learning and exploration  via the notion of {\it randomized policies}, known  in the control literature as  {\it relaxed controls} and in the game theory as {\it mixed strategies}. These are policies, say $\pi_t$, of  the representative agent with  $\pi_t\in \PP(U)$, where the action space $U$ is a closed subset of a Euclidean space and  $\PP(U)$ 
is the set of density functions of probability measures on $U$ that are absolutely continuous with respect to the Lebesgue measure.
 Namely, $\pi_t \in \PP(U)$ if and only if
\begin{eqnarray}\label{PU}
\int_U \pi_t(u)du=1\,\,\,\,\,\, \textit{ and }\,\,\,\,\,\, \pi_t(u)\ge 0  \textit{ a.e. on } U.
\end{eqnarray}

The third ingredient is the
{\it entropy regularization}, which is adopted to encourage 
exploration.  For this, we will use both the Shannon entropy 
and the cross-entropy, denoted by $\HH_{SE}$ and  $\HH_{CE}$, respectively (see \eqref{eq:shannon} and \eqref{eq:cross}).

\paragraph{Controlled state process with randomized policies.}

We incorporate the above components in a finite horizon setting $[0,T]$, $0<T<\infty$. For this, we introduce $\mu:=\{\mu_s,\,\, s\in[t,T]\}$ to be the flow of population state distribution with $\mu_s \in \mathcal{P}(\mathbb{R})$ and $\alpha:=\{\alpha_s,\,\,s\in[t,T]\}$ to be the flow of population action distribution with $\alpha_s:\mathbb{R}\rightarrow\mathcal{P}(U)$, starting from time $t\in[0,T]$. $\alpha_s(\cdot;x) \in \mathcal{P}(U)$ represents the action distribution of the population at state $x\in \mathbb{R}$.
 Occasionally, $\alpha$ and $\mu$ will  be also  called the {\it mean field information}.

Next we define the controlled state process of the representative agent. Given $t\in[0,T]$ and 
exogenous flows, say $\alpha$ and  $\mu$ with $\mu_t = \nu$, the representative agent adopts a randomized policy $\pi=\{\pi_s \in \mathcal{P}(U), s \in [t,T]\}$  over an admissible policy set ${\cal A}$ (to be specified below). Then, following the paradigm recently proposed in \cite{WZZ2018}, the controlled state process is assumed to follow
 \begin{eqnarray}\label{eq:state}
d X^{\pi}_s  &&= \left(\int_U b(s,X^{\pi}_s,\mu_s,\alpha_s,u)\pi_s(u)du \right)ds + \left(\sqrt{\int_U \sigma^2(s,X^{\pi}_s,\mu_s,\alpha_s,u)\pi_s(u)du }\right)dW_s,  \\
  X^{\pi}_t && = \xi \sim \nu, \ \ \mu_t = \nu, s \in [t, T]. \nonumber
\end{eqnarray}
Here  ${W}=\{W_t\}_{t\in[0,T]} $ is a standard  Brownian motion defined on a  filtered probability space $(\Omega,\mathcal{F},\{\mathcal{F}_t\}_{ t\in [0,T]},\P)$, with $\{\mathcal{F}_t\}_{t\in[0,T]}$  satisfying the usual conditions; $\nu \in \mathcal{P}(U)$  is the distribution of the initial state satisfying $\int x^2\nu(dx)<\infty$; $\xi$ is  a random variable independent of $W$ and $\mathcal{F}_t$-measurable; and $b, \sigma: [0,T]\times \mathbb{R}\times \mathcal{P}(\mathbb{R}) \times \mathcal{P}(U) \times U \hookrightarrow \mathbb{R}$. 

We note the particular form of the state process \eqref{eq:state} is a consequence of the  aggregation of $\widehat{X}_s^u$ over action $u_s\in U$ where
\begin{eqnarray*}
d \widehat{X}^u_s = {b}\left(s,\widehat{{X}}^u_s,\mu_s,\alpha_s,u_s\right) ds +{\sigma}\left(s,\widehat{{X}}^u_s,\mu_s,\alpha_s,u_s\right) dW_s.
\end{eqnarray*}
Such policies  $u_s\in U$ are also called {\it pure strategies} in game theory.
Pure strategies  and mixed strategies are closely related. Indeed, $u = \{u_s, s \in [0,T]\}$ can be regarded as a Dirac measure $\pi = \{\pi_s(u), s \in [0,T]\}$ where $\pi_s(\cdot) = \delta_{u_s}(\cdot)$.  In this case $\pi_s$ does not have a density, and hence $\pi_s\notin \mathcal{P}(U)$. (We refer the readers to \cite{WZZ2018} for more details).
 
\paragraph{Game payoff with entropy regularization.} 
The objective of the representative agent is to maximize her payoff function $J$ and solve for
\begin{equation}\label{mfg_value_general}
   V(t \,|\,\mu,\alpha)= \sup_{\pi\in \mathcal{A}} J(t,\pi\,|\,\mu,\alpha),  
 \end{equation}
where the     entropy-regularized payoff is defined as
 \begin{eqnarray} \label{mfg_obj_general}
&& 
J(t,\pi\,|\,\mu,\alpha)\nonumber=\E\left[ \int_t^{T} \left(\, \int_U \Big(r(s,X^{\pi}_s,\mu_s,\alpha_s,u)\pi_s (u) du+\lambda_{SE} \mathcal{H}_{SE}(\pi_s)   +\lambda_{CE} \mathcal{H}_{CE}(\pi_s,\alpha_s,\mu_s)\Big) \right)ds\right. \nonumber\\
\vspace{15 mm}\nonumber
\\
&&\left.\qquad \qquad \qquad +g(X^{\pi}_T,\mu_T,\alpha_T)\Big\vert  \mu,\alpha\right].\label{mfg}
\end{eqnarray}
The Shannon entropy $\HH_{SE}$ and cross-entropy $\HH_{CE}$ are defined  as 
\begin{eqnarray}
\mathcal{H}_{SE}(\pi_s) &=& {-}\int_U\pi_s (u)\ln \pi_s (u)du, \,\,\pi_s \in \mathcal{P}(U),\label{eq:shannon}\\
\vspace{15 mm}\nonumber
\\
\mathcal{H}_{CE}(\pi_s,\alpha_s,\mu_s) &=&{-} \int_U\pi_s (u)\int \ln \alpha_s (u;x)\,\mu_s(dx)du,\,\,\,\pi_s \in \mathcal{P}(U).\label{eq:cross}
\end{eqnarray}
In addition, $r: [0,T]\times \mathbb{R}\times \mathcal{P}(\mathbb{R})\times \mathcal{P}(U) \times U \hookrightarrow \mathbb{R}$ and $g:  \mathbb{R}\times \mathcal{P}(U)\times \mathcal{P}(U) \hookrightarrow \mathbb{R}$ are the running reward and terminal reward functions of the representative agent, while  $\lambda_{SE} >0$  is the (temperature) parameter to control the degree of self-exploration  and $\lambda_{CE} \ge 0 $  is the (temperature) parameter to control the degree of exploration on  the actions of the other agents. From an information-theoretic perspective, $\lambda_{SE}\mathcal{H}_{SE}$ and $\lambda_{CE}\mathcal{H}_{CE}$ quantify the information gain from exploring the unknown environment and the policies chosen by the other agents.

\paragraph{Observable quantities.} In a game with learning, the functions
$b$, $\sigma$, $r$ and $g$ are {\it unknown}. The representative agent takes actions while interacting with (the continuum) of the other agents. This interaction takes several rounds. 

In each round starting from time $0$, the agent observes $\{\alpha_s\}_{s\in[0,t]}$, $\{\mu_s\}_{s\in[0,t]}$ and $\{X^{\pi}_s\}_{s\in[0,t]}$ at time $t\in [0,T]$; the reward will not be revealed until time $T$, the end of each round;  at time $T$, she will observe the {\it realized} cumulative reward $\widehat{j} \left(0,\pi|\alpha,\mu\right)$ with
\begin{eqnarray*}
\widehat{j} \left(0,\pi|\alpha,\mu\right) &:=&\int_0^{T} \Big( \int_U \left(r(s,X^{\pi}_s,\mu_s,\alpha_s,u)\pi_s (u) du+\lambda_{SE} \mathcal{H}_{SE}(\pi_s)   +\lambda_{CE} \mathcal{H}_{CE}(\pi_s,\alpha_s,\mu_s) \Big)\right)ds\\ 
\vspace{15 mm}\nonumber
\\
&&+g(X^{\pi}_T,\mu_T,\alpha_T),
\end{eqnarray*}
which is associated with the corresponding {\it single} trajectory $\{X^{\pi}_s\}_{s\in[0,T]}$ under policy $\pi$ and the population behavior $\{\alpha_s\}_{s\in[0,T]}$, $\{\mu_s\}_{s\in[0,T]}$ in this round. Note that $\widehat{j} \left(0,\pi|\alpha,\mu\right)$ is one realized sample reward, which is different from the expected reward in \eqref{mfg_obj_general}.

\paragraph{Admissible policies.}  A policy $\pi  \in \mathcal{A}(t,\mu,\alpha)$ is admissible if
\begin{enumerate}
\item[(i)] for each $s \in [t,T]$, $\pi_s \in \mathcal{P}(U)$ a.s.;
\item[(ii)] for each ${Z} \in \mathcal{B}(U)$
with $\mathcal{B}(U)$ being the Borel algebra on $U$,
$\{\int_{Z} \pi_s(u) du, s\in [t,T]\}$ is $\mathcal{F}_t$-progressively measurable;
\item[(iii)] the SDE  \eqref{eq:state} admits a unique strong solution $X^{\pi}:=\{X^{\pi}_s,s \in [t,T]\}$,
with $\pi$ being used; 
\item[(iv)] the expectation on the right hand side of \eqref{mfg} is finite;
\item[(v)] there exists a measurable function ${\widetilde{\pi}}: [t, T ]\times
\R \rightarrow \PP(U)$ such that $$\P\Big(\pi_s(du) = {\widetilde{\pi}_s}(du;X^{\pi}_s), \quad \forall s\in [t,T] \Big) = 1.$$
\end{enumerate}

Condition (v) imposes that the admissible policy is Markovian, i.e., closed-loop policy in feedback form.
\paragraph{Alternative formulation of the MFG with learning.}
We note  that  problem \eqref{mfg_value_general} treats the initial state $\xi$ as a genuine source of randomness, in addition to the stochasticity  from the Brownian motion $W$. Frequently, the following alternative interpretation, with a {\it deterministic} initial state $x$  is useful for solving analytically the MFG. Specifically, let
\begin{eqnarray} \label{mfg_obj_alternative}
\widetilde{V}(t,x\,|\,\mu,\alpha):=&& \sup_{\pi\in {\cal A}}\widetilde{J}(t,x\,|\,\pi,\mu,\alpha)\nonumber\\
:= &&\E\left[ \int_t^{T} \left( \int_U \left(r(X^{\pi}_s,\mu_s,\alpha_s,u)\pi_s (u) du+\lambda_{SE} \mathcal{H}_{SE}(\pi_s)   +\lambda_{CE} \mathcal{H}_{CE}(\pi_s,\alpha_s,\mu_s) \right)\right)ds\right. \nonumber\\
\vspace{15 mm}\nonumber
\\
&&\left.\ \ \ \ \ \ \ \ \ \ \ \ \ +g(X^{\pi}_T,\mu_T,\alpha_T)\Big\vert X^{\pi}_t = x, \mu,\alpha\right],
\end{eqnarray}
subject to 
\begin{eqnarray}\label{eq:mfg_dynamics}
d X^{\pi}_s  &&= \left(\int_U b(s,X^{\pi}_s,\mu_s,\alpha_s,u)\pi_s(u)du \right)ds + \left(\sqrt{\int_U \sigma^2(s,X^{\pi}_s,\mu_s,\alpha_s,u)\pi_s(u)du }\right)dW_s,  \\
 \vspace{15 mm}\nonumber
\\
  X^{\pi}_t && = x, \ \ \mu_t = \nu, \,\,s \in [t, T]. \nonumber
\end{eqnarray}
Then, it easily follows that
$$\E_{\xi\sim\nu}[\widetilde{V}(t,\xi\,|\,\mu,\alpha)] =  V(t\,|\,\mu,\alpha).$$
While conceptually this approach is less general, it is frequently used - as in \cite{LZ2017} and herein - to solve the MFG explicitly.


\paragraph{Nash Equilibrium (NE) for MFG with learning.}
To analyze game \eqref{eq:state}-\eqref{mfg_value_general}, we adopt the well-known NE criterion. 

\begin{definition}[NE for MFG] 
\label{def:NE} For game \eqref{mfg_value_general} with an initial state distribution $\nu$ and state process \eqref{eq:state}, 
an agent-population profile $(\pi^{*},\mu^*,\alpha^*):=\{(\pi^{*}_s,\mu^*_s,\alpha^*_s), t\le s \le T\}$ is called NE if the following conditions hold:
\begin{enumerate}
\item[{\bf A.}] (Single-agent-side) For the fixed population state-action distribution $(\mu^*,\alpha^*)$ and any policy $\pi\in \mathcal{A}$,
$$J(t,\pi\,|\,\mu^{*},\alpha^*) \leq J(t,\pi^{*}\,|\,\mu^*,\alpha^*).$$
\item[{\bf B.}] (Population-side)  $\pi_s^{*}(u;x) = \alpha_s^*(u;x)$, for all $x\in \mathbb{R}$. In addition, $\P_{X^*_s}=\mu_s^*$ for any $s \in [t,T]$, where $X^*$ solves \eqref{eq:state} when policy $\pi^*$ is adopted with the initial population state distribution $\mu_t^* = \nu$.
\end{enumerate}
Given a NE $(\pi^{*},\mu^*,\alpha^*)$,  $$ V(t\,|\,
\mu^*,\alpha^*):=J(t,\pi^*\,|\,\mu^{*},\alpha^*) = \max_{\pi\in \AA}J(t,\pi\,|\,\mu,\alpha^*)$$
is called  a {\it game value} associated with this NE.
\end{definition}
Given  $(\mu^*,\alpha^*)$,  condition {\bf A} captures the optimality of $\pi^{*}$ while  condition {\bf B}  ensures the consistency of the solution so that  the state and action flows of the single agent  match those of the population. Note that   uniqueness of NE for MFG is, in general, rare when mixed strategies are allowed (see, for example, \cite{Lacker2015}).

{
\paragraph{Solvability of MFG.}
There are three classical approaches to show the existence of MFG solution with pure strategies (or strict controls): the PDE  (i.e., three-step fixed point) approach \cite{LL2007,HCM2007}, the probabilistic approach \cite{carmona2013mean,carmona2015probabilistic} and the master equation approach \cite{cardaliaguet2019master}. The uniqueness of the MFG solution with pure strategies can be verified with certain technical conditions such as the small parameter conditions \cite{LL2007} or the monotonicity condition \cite{HCM2007}. 

In the framework of MFG with relaxed controls, we follow the    three-step fixed point approach to solve \eqref{eq:state}-\eqref{mfg_value_general}: 
 \begin{itemize}
\item {\bf Step 1:} Fix a population state-action distribution $(\mu,\alpha)$ and an initial state $x$. Then, solving the MFG \eqref{eq:state}-\eqref{mfg_value_general} is reduced to solving a stochastic control problem with randomized policies (relaxed controls).
\item {\bf Step 2:} Let $X_s^{\pi,x}$ be the controlled state process under the optimal  policy $\pi$ from the initial state $x$ in Step 1. Update $\alpha^{\prime}_s(\cdot,y) =\pi_s(\cdot,y) $ for all $y \in \R$ and $s \in [t,T]$.  Denote $X_s^{\pi,\xi}$ the  controlled state process under $\pi$ from some random initial state $\xi\sim \nu$. Then, update $\mu^{\prime}_s=\mathbb{P}_{X_s^{\xi,\pi}}$. 
\item {\bf Step 3:} Repeat Steps 1 and 2 until $(\mu^{\prime},\alpha^{\prime})$ converges.
\end{itemize}
\smallskip

Note that there is no guarantee that the above procedure will yield any MFG solution since {\bf Step 1} may have multiple solutions under relaxed controls. Moreover, by the nature of relaxed controls, the candidate fixed point(s) would be the fixed point(s) of a set-valued map as described in \cite{Lacker2015}.  Nevertheless, for a family of linear-quadratic MFG, which will be introduced in Section \ref{sec:LQ}, 
one can build proper verification arguments to show that the explicit fixed-point solution is indeed a solution to the MFG problem \eqref{eq:state}-\eqref{mfg_value_general}.

In general, the uniqueness of the MFG solution with relaxed controls does not hold unless there are  additional convexity properties of the value function (see, for example, \cite{DL}). Here, the convexity in the linear-quadratic framework fails to hold when entropy regularization is included.
}

\section{Shannon Entropy and Enhanced Entropy for MFG with Learning}
\label{sec:LQ}
In the mathematical formulation for MFG with learning of Section \ref{sec:setup}, we 
 analyze the information theoretic gain for  two types of entropies: Shannon entropy $\mathcal{H}_{SE}$ and enhanced entropy, which is a linear combination of Shannon entropy  and cross-entropy $\lambda_{SE}\mathcal{H}_{SE}+\lambda_{CE}\mathcal{H}_{CE}$, with temperature parameters $\lambda_{SE}$ and $\lambda_{CE}$.
 We  study the impact of this entropy regularization {within a class of LQ games} in a finite time horizon. 
 LQ games are the building blocks of stochastic games and often bring critical  insights from their closed-form solutions
 (\cite{bardi2012explicit,
 bensoussan2016linear}). 
Among others, we will see that the LQ games we analyze yield  time-dependent optimal  policies, with  
 time-dependent Gaussian efficient explorations.

  \subsection{Game with Shannon Entropy}
\label{subsec:entropy}
We start with the case of using only Shannon entropy for exploration, namely

\begin{eqnarray}
\label{MFG-SE}
V_{SE}(t \,|\,\mu):=&&  \sup_{\pi\in {\cal A}}  \ \ \ J_{SE}(t,\pi \,|\,\mu)   \nonumber \\
:=&&\sup_{\pi\in {\cal A}} \E \left[ \int_t^{T}  \left( \int_{\mathbb{R}} -\frac{Q}{2}(X^{\pi}_s-m_{s})^2\pi_s (u) du+\lambda_{SE} \HH_{SE}(\pi_s)
\right.\right)ds  \left.\left.-\frac{\bar{Q}}{2}(X^{\pi}_T-m_T)^2 
 \right\vert \mu\right], \nonumber
 \end{eqnarray}
subject to
\[d X^{\pi}_s = \left(\int_{\mathbb{R}} (A(m_s-X^{\pi}_s)+Bu)\pi_s(u)du \right)ds + D\left(\sqrt{\int_{\mathbb{R}} u^2 \pi_s(u)du}\right) dW_s,  \ \ 
X^{\pi}_t=\xi\sim \nu. \tag{\textbf{MFG-SE}}
\]
Here $\mu_t = \nu$, and $m_s=\int x \mu_s(dx)$  ($s \in [t,T]$).  We assume $A>0$, $Q>0$, $\bar{Q}>0$, and $\lambda_{SE}>0$. We take the action space to be $U=\mathbb{R}$, and without loss of generality,  $B>0$ and $D>0$.

We remark  that $\alpha := \{\alpha_s\}_{s\in[t,T]}$ 
does not appear  in the game formulation ({\bf{MFG-SE}}). This is because when  only Shannon entropy is incorporated,  there is no interaction between the policy of the representative agent  and the  population action distribution $\alpha$.

There are  two types of rewards in this game: the running reward $-\frac{Q}{2}(X^{\pi}_s-m_s)^2$  that penalizes any deviation from the current average state of the population at time $s\in[t,T]$, and the terminal reward $-\frac{\bar{Q}}{2}(X^{\pi}_T-m_T)^2$ that penalizes deviation from the  average state of the population at terminal time $T$.  There are also two types of interaction:   the real time interaction $A(m_s-X^{\pi}_s)$ and $\frac{Q}{2}(X^{\pi}_s-m_s)^2$ for $s\in[t,T]$, and the interaction at terminal time  $\frac{\bar{Q}}{2}(X^{\pi}_T-m_T)^2$.\\

Next, we present one of the main results herein which provides an explicit NE solution for the MFG we consider.  For notational convenience, We denote by $\mathcal{N} (\cdot|\nu, \sigma^2)$ the density function of a
Gaussian random variable with mean $\nu$ and variance $\sigma^2$.
\begin{theorem}[{\bf MFG-SE}]\label{thm:current_shannon}  Let $m^* =\mathbb{E}[\xi]$ and 
\begin{equation}
\label{eqn:V}
\widetilde{V}_{SE}(t,x)=-\frac{\eta^{SE}_t}{2}(x-m^*)^2+\gamma^{SE}_t,
\end{equation} with 
\begin{equation}
\label{eqn:eta}
\eta^{SE}_t = {\bar{Q}}\exp\left(-\left(2A+\frac{B^2}{D^2}\right)(T-t)\right)+{\frac{Q}{2A+\frac{B^2}{D^2}}\left(1- \exp\left({-\left(2A+\frac{B^2}{D^2}\right)(T-t)}\right)\right)}>0,
\end{equation}
and $$\gamma^{SE}_t =  \frac{\lambda_{SE}}{2}\ln\left({\frac{2\pi  \lambda_{SE}}{D^2}}\right)(T-t)-\int_t^T\frac{\lambda_{SE}}{2}\ln\left(\eta^{SE}_z\right)dz.$$
Then, $$V^*_{SE}(t) := \E_{\xi\sim\nu}[\widetilde{V}_{SE}(t,\xi)]$$ is a game value of {\bf (MFG-SE)} associated with the NE policy
\begin{equation}\label{eqn:pi}
\pi_{s}^{SE*}(u;x) = \mathcal{N}\left(u\left|\frac{B(m^*-x)}{D^2},\frac{\lambda_{SE}}{D^2\eta^{SE}_s}\right)\right., \ \  s\in [t,T].
\end{equation}

The corresponding controlled state process under \eqref{eqn:pi} is the unique solution of the SDE, 
\begin{eqnarray}\label{eq:se_optimal_dynamics}
d X_s^* &=& \left(A+\frac{B^2}{D^2}\right)(m^*-X_s^*)ds+\sqrt{\frac{B^2(X_s^*-m^*)^2}{D^2}+ \frac{\lambda_{SE}}{\eta^{SE}_s}}\,dW_s, \\
\vspace{15 mm}\nonumber
\\
 X_t^{*} &=&\xi \sim \nu, \,\,s\in[t,T]. \nonumber
\end{eqnarray}
 In addition,  the mean state of the population under policy \eqref{eqn:pi} is time-independent, i.e., 
\begin{equation}\label{eqn:m}
m_s^* = \mathbb{E}[X_s^*] = m^*,  \  \  s \in [t,T].
\end{equation}
\end{theorem}

\begin{remark}
Theorem \ref{thm:current_shannon}  provides important guidance for exploration from an information-theoretic perspective. 
  It suggests that, with Shannon entropy regularization, the associated optimal policy  $\pi_{s}^{SE*}(u;x)$  from \eqref{eqn:pi}  is Gaussian,  mean-reverting and with time-dependent  variance. This is  useful for MARL algorithmic design as  the  agent can now focus on a much smaller class of  policies 
\begin{eqnarray}\label{eqn:hatcontrol}
\widehat{\pi}_s(u;x) \sim \mathcal{N}\left(\hatM(m_s-x),\hatV_s\right),
\end{eqnarray}
with $m_s = \int_{\mathbb{R}} x \mu_s(dx)$,  $\widehat{M}$ some scalar and  $\widehat{\sigma}^2 =\{\widehat{\sigma}^2_s\}_{s\in[t,T]}$ a variance exploration process. Meanwhile, she  can improve her estimate on $\hatM$ and $\widehat{\sigma}^2$ of the above policy while interacting with the system and other agents, and observing the outcome at the end of each round of play. 
Indeed, notice that the controlled state process  becomes 
\begin{eqnarray}\label{eqn:hatdynamics}
d X^{\widehat{\pi}}_s = \left(A+B\hatM\right)(m_s-X^{\widehat{\pi}}_s)ds + \left(D\sqrt{\hatM^2(X^{\widehat{\pi}}_s-m_s)^2 +\hatV_s} \right)d W_s,
\end{eqnarray}
with $X_t =\xi$.  Thus,  the following simple corollary will be useful for MARL (see also more details in Section \ref{sec:experiment} where this result is used for algorithm design).
\begin{corollary}\label{cor:random_policy}
If the representative agent follows policy  \eqref{eqn:hatcontrol} under a given mean field information $\mu=\{\mu_s\}_{s\in [t,T]}$, then the payoff is given by
\begin{eqnarray}\label{hat_cost}
J_{SE}(t,\widehat{\pi}\,|\,\mu) &=&-\frac{Q}{2} \int_t^T \left(\phi_s^2 -2{m}_s\,\widehat{m}_s\,+{m}^2_s\right)ds\nonumber\\
&& +\frac{\lambda_{SE}}{2} \int_t^T \ln \left(2\pi e \hatV_s\right)ds - \frac{\bar{Q}}{2} \left(\phi_T^2 -2{m}_T\,\widehat{m}_T\,+{m}^2_T\right),
\end{eqnarray}
where
\[
\phi_s^2 = e^{\left(2\hatK +D^2\hatM^2\right)(s-t)}\left(\E[\xi^2] + \int_t^s e^{- D^2\hatM^2(z-t)}d(z)dz\right),
\]
with
\[
d(s) = -2\E[\xi]\, e^{- \hatK (s-t)}\hatK m_s+\left(\int_t^s e^{- \hatK (z-t)}\hatK m_zdz\right) e^{- \hatK (s-t)}\hatK m_s +e^{-2 \hatK (s-t)}D^2\left(\hatM^2m_s^2-2\hatM^2m_s \widehat{m}_s +\hatV_s \right),
\]
\begin{eqnarray*}
m_s = \int_{\mathbb{R}} x\mu_s(dx),\quad \widehat{m}_s &=&e^{\hatK (s-t)}\E [\xi] + \int_t^s e^{\hatK (s-z)}\hatK m_z\,dz,\quad \mbox{and} \ \ \hatK = -\left(A+B\hatM\right).
\end{eqnarray*}
\end{corollary} 

\end{remark}

Next, we analyze the game with an additional cross-entropy regularization.
\subsection{Game with Enhanced Entropy (Linear Combination of Shannon Entropy and Cross-entropy)}
  The objective of this game is  to find
\begin{eqnarray*}
V_{EE}(t \,|\, \mu, \alpha):=&& \sup_{\pi\in {\cal A}}J_{{EE}}(t, \pi\,|\,\mu,\alpha) \\
:= &&  \sup_{\pi\in {\cal A}}\E \left[ \int_t^{T} \left(  -\frac{Q}{2}(X^{\pi}_s-m_{s})^2+\lambda_{SE} \HH_{SE}(\pi_s)+\lambda_{CE} \HH_{CE}(\pi_s,\alpha_s,\mu_s)\right)ds\right.\\
\vspace{15 mm}\\
&&\qquad \qquad\left.\left.-\frac{\bar{Q}}{2}(X^{\pi}_T-m_T)^2  \right\vert\,\alpha,\mu\right],
\end{eqnarray*}
subject to
\[d X^{\pi}_s = \left(\int_\mathbb{R} (A(m_s-X^{\pi}_s)+Bu)\pi_s(u)du \right)ds + D\sqrt{\int_\mathbb{R} u^2\pi_s(u)du} dW_s, \ \ X^{\pi}_t = \xi\sim\nu.  \tag{\textbf{MFG-EE}}\]
Here $\mu_t = \nu$, $m_s=\int x \mu_s(dx)$, $s\in [t,T]$, and we assume $\lambda_{SE}>0$, $\lambda_{CE} \ge 0$, $Q>0$, $\bar{Q}>0$, and $A>0$. Without loss of generality, we also take $B>0$ and $D>0$.

\begin{theorem}[{\bf MFG-EE}]\label{thm:current_cross}  Let $m^*=\E[\xi]$, $\rm{Var}[\xi] = \mathbb{E}[\xi^2] -(\mathbb{E}[\xi])^2$, and 
\begin{equation}\label{eqn:V2}
\widetilde{V}_{EE}(t,x):=-\frac{\eta^{EE}_t}{2}(x-m^*)^2+\gamma^{EE}_t,
\end{equation} 
with
\begin{eqnarray}\label{eqn:eta2}
\eta^{EE}_t &=& {\bar{Q}}\exp\left(-\left(2A+\frac{B^2}{D^2}\frac{\lambda_{SE}+\lambda_{CE}}{\lambda_{SE}}\right)(T-t)\right)\nonumber\\
\vspace{15 mm}\nonumber
\\
&&\quad+\frac{Q}{2A+\frac{B^2}{D^2}\frac{\lambda_{SE}+\lambda_{CE}}{\lambda_{SE}}}\left(1-\exp\left(-\left(2A+\frac{B^2}{D^2}\frac{\lambda_{SE}+\lambda_{CE}}{\lambda_{SE}}\right)(T-t)\right)\right),
\end{eqnarray}
{where $\eta^{EE}_t>0$, for $t\in[0,T]$,}
and
\begin{eqnarray*}
\gamma_t^{EE} &=&  \frac{\lambda_{SE}+\lambda_{CE}}{2}\ln\left({\frac{2\pi ( \lambda_{SE}+\lambda_{CE})}{D^2}}\right)(T-t)-\frac{\lambda_{SE}+\lambda_{CE}}{2}\int_t^T\ln(\eta^{EE}_z)dz\\
\vspace{15 mm}\\
&& +\frac{B^2}{2D^2}\frac{\lambda_{CE}(\lambda_{SE}+\lambda_{CE})}{\lambda_{SE}^2}\int_t^T\eta^{EE}_z\,\kappa^{EE}_z dz,
\end{eqnarray*}
 with
 \begin{eqnarray}
\kappa^{EE}_s=e^{ (2K+M)(s-t)} \rm{Var}[\xi] + \int_t^s e^{(M+2K) (s-z)}\frac{\lambda_{SE}+\lambda_{CE}}{\eta^{EE}_z} dz
 \end{eqnarray}
and
\begin{eqnarray*}
 K =-\left( A+\frac{B^2}{D^2}\right){\frac{\lambda_{SE}+\lambda_{CE}}{\lambda_{SE}}},\quad M = \left(\frac{B}{D}{\frac{\lambda_{SE}+\lambda_{CE}}{\lambda_{SE}}}\right)^2.
\end{eqnarray*}
Then,
  $$V^*_{EE}(t):=\E_{\xi\sim\nu}\left[\widetilde{V}_{EE}(t,\xi\,|\,\mu^*,\alpha^*)\right]$$ is a game value of {\bf (MFG-EE)}, with the associated NE policy 
\begin{eqnarray}\label{eqn:pi2}
\pi_{s}^{EE*}(u;x) = \mathcal{N}\left(\frac{\lambda_{SE}+\lambda_{CE}}{\lambda_{SE}}\frac{B(m^*-x)}{D^2},\frac{\lambda_{SE}+\lambda_{CE}}{D^2\eta^{EE}_s}\right).
\end{eqnarray}
Furthermore, the optimal controlled state process $X_s^*$ under policy \eqref{eqn:pi2} is the unique solution of the SDE, 
\begin{eqnarray}
d X_s^* &&= \left(A+\frac{\lambda_{SE}+\lambda_{CE}}{\lambda_{SE}}\frac{B^2}{D^2}\right)(m^*-X_s^*)ds\nonumber\\
\vspace{15 mm}\nonumber\\
&&+D\sqrt{\left(A+\frac{\lambda_{SE}+\lambda_{CE}}{\lambda_{SE}}\frac{B^2}{D^2}\right)^2(X_s^*-m)^2+ \frac{\lambda_{SE}+\lambda_{CE}}{D^2\eta^{EE}_s}}dW_s, \label{eq:ee_optimal_dynamics}\\
\vspace{15 mm}\nonumber
\\
X_t^* && = \xi\sim\nu, \,\,s\in [t,T].\nonumber
\end{eqnarray}
 In addition, $\mu_t^* = \mathbb{P}_{X_t^*}$, $\alpha^*_s(u;x) =\pi_{s}^{EE*}(u;x) $, and the mean state of the population under policy \eqref{eqn:pi2} is time independent, i.e., 
\begin{equation}\label{eqn:m2}
m_s^* = \mathbb{E}[X_s^*] = m^*,  \  \  s \in [t,T].
\end{equation}
\end{theorem}
\smallskip

Before providing the proof, a few remarks are in place.
\subsection{Discussion.} 
\label{sec:discussion}
In both linear-quadratic MFGs, with either only the Shannon entropy ({\bf MFG-SE}) or with the  additional  cross-entropy ({\bf MFG-EE}), there are several similarities.
\begin{itemize}
\item The form of the optimal policies \eqref{eqn:pi} and \eqref{eqn:pi2} suggests that Gaussian exploration  is optimal  when  entropy regularization is  introduced in the MFG with learning.
This is consistent with recent  works of \cite{WZZ2018, WZ2019} for continuous-time single-agent RL and is also supported by the empirical studies of \cite{LHPHTSW2015} and \cite{plappert2017}.  
\item 
Both the means of the optimal policy $\pi^{SE*}_s(u;x)$  in \eqref{eqn:pi}  and the  optimal policy $\pi^{EE*}_s(u;x)$ in \eqref{eqn:pi2}  are  influenced by both the mean field interaction and the current state of the representative agent.  On the other hand, both their  variances   are time-dependent.

In addition,  the strength of their mean reversion is quantified by the coefficient $\frac{B}{D^2}$, which indicates that a smaller variance signifies less uncertainty in the game, hence a faster mean reverting policy.
\item Equation \eqref{eqn:eta} for $\eta_s^{SE}$ and equation \eqref{eqn:eta2} for $\eta_s^{CE}$ suggest that when time $s$ is sufficiently small,  the term 
${\frac{Q}{2A+\frac{B^2}{D^2}}\left(1- \exp\left({-\left(2A+\frac{B^2}{D^2}\right)(T-s)}\right)\right)}$  dominates  $\eta^{SE}_s$, whereas  
$\eta^{EE}_s$ is dominated by
$\frac{Q}{2A+\frac{B^2}{D^2}\frac{\lambda_{SE}+\lambda_{CE}}{\lambda_{SE}}}\left(1-\exp\left(-\left(2A+\frac{B^2}{D^2}\frac{\lambda_{SE}+\lambda_{CE}}{\lambda_{SE}}\right)(T-s)\right)\right)$. Thus, when time $s$ is small, the cost of exploration is low and the representative agent has more incentive to explore in upcoming times.

Conversely, when time $s$ is sufficiently large and especially when $s\sim T$,  $\eta^{SE}_s$ is dominated by  the term
 ${\bar{Q}}\exp\left(-\left(2A+\frac{B^2}{D^2}\right)(T-s)\right)$, whereas
 ${\bar{Q}}\exp\left(-\left(2A+\frac{B^2}{D^2}\frac{\lambda_{SE}+\lambda_{CE}}{\lambda_{SE}}\right)(T-s)\right)$ dominates
 $\eta^{EE}_s$. 
 Thus, the cost of exploration increases  as time  $s$ approaches $T$.  This implies  that the agent is more sensitive to the terminal reward  and explores less when the game approaches termination.
 \item In the very special case $A\equiv Q\equiv 0$,  there is no intermediate payoff. Then, the variance of $\pi_s^{SE*}$ and $\pi_s^{EE*}$ decreases when time $s$ increases, implying  more exploration at the very  beginning and less towards the very end.
\end{itemize}
\smallskip

Despite the above similarities, there is an important difference: 
\begin{itemize}
\item  The Shannon entropy  and the cross-entropy affect the optimal policy $\pi_s^{SE*}$   and $\pi_s^{EE*}$ differently.
Indeed,  the mean of the optimal policy  $\pi_s^{EE*}(u;x)$  depends on the ratio between $\lambda_{CE}$ and $\lambda_{SE}$, while $\lambda_{SE}$ and $\lambda_{CE}$ impact the variance of $\pi_s^{EE*}(u;x)$  through both the $\frac{\lambda_{CE}}{\lambda_{SE}}$ and $\lambda_{SE}+\lambda_{CE}$ terms.
 In particular, with the additional cross-entropy, one will  explore more and, consequently, the learning  procedure would converge faster.

\end{itemize}

\subsection{Derivations and Proofs of Main Results}
The solution approach consists of two steps. The first is to find a candidate  solution based on the classical fixed-point approach introduced in Section \ref{sec:setup}. The second  is to verify the candidate solution via a verification theorem. 

\paragraph{NE Derivation of {\bf (MFG-SE)}.} To ease the exposition, we drop the subscript SE.
\begin{proof}{Proof of Theorem \ref{thm:current_shannon}.} For a given {admissible} policy $\pi \in \mathcal{A}$, the forward equation for $p(s,x)$, the density of $X_s$, is given by,
\[
\partial_s p(s,x) = -\partial_x \left(\left(A(m_s-x)+B \int_\mathbb{R} u\, \pi_s(u;x)du \right)p(s,x)\right) +\frac{1}{2}\partial_{xx}\left(p(s,x)\int_\mathbb{R} D^2u^2 \,\pi_s(u;x)du \right),
\]
with initial density $p(t,x) = \nu(x)$. Here, $m_s = \int x p(s,x)dx$, $s \in [t,T]$.  

We first proceed heuristically with the associated  HJB equation, derive a solution, and then  validate this solution through a verification argument.

\noindent {\bf Step 1 (solving the control problem):} Given fixed mean-field information $\{m_s\}_{ s \in [0,T]}$ which is deterministic, the HJB equation for the value function $\widetilde{V}(s,x)$ can be written as
\begin{eqnarray}\label{hjb}
-\partial_s \widetilde{V}(s, x)=
&&\max_{\pi_s \in \PP(\mathbb{R})} \left(\left(A(m_s-x)+B \int_\mathbb{R} u \pi (u;s,x)du\right) \widetilde{V}_x(s,x) \right.\\
\vspace{15 mm}\nonumber
\\
&&\,\,\,\,\left.-\frac{Q}{2}(m_s-x)^2-\lambda_{SE} \int_\mathbb{R} \pi_s(u;x)\ln \pi_s(u;x)du +\frac{1}{2}\left(\int_\mathbb{R} D^2u^2 \pi_s(u;x)du\right)\partial_{xx}\widetilde{V}(s,x) \right),\nonumber
\end{eqnarray}
with terminal condition $\widetilde{V}(T,x) = -\frac{\overline{Q}}{2}(x-m_T)^2$.
Recall that $\pi_s(u;x) \in \mathcal{P}(U)$ if and only if \eqref{PU} holds. Solving the constrained maximization problem on the right hand side of \eqref{hjb} yields
\begin{eqnarray*}
\pi_s^*(u;x)= \frac{\exp\left(\frac{1}{\lambda_{SE}}\left(-\frac{Q}{2}(x-m_s)^2+\frac{1}{2}D^2u^2\partial_{xx}\widetilde{V}+(A(m_s-x)+Bu)\widetilde{V}_x\right)\right)}{\int_\mathbb{R} \exp\left(\frac{1}{\lambda_{SE}}\left(-\frac{Q}{2}(x-m_s)^2+\frac{1}{2}D^2u^2\partial_{xx}\widetilde{V}+(A(m_s-x)+Bu)\widetilde{V}_x\right)\right)du}.
\end{eqnarray*}
Thus, the optimal  policy  is expected to be \textit{Gaussian} with mean
$\frac{B \partial_x\widetilde{V}}{-D^2\partial_{xx}\widetilde{V}}$ and variance $\frac{\lambda_{SE}}{-D^2\partial_{xx}\widetilde{V}}$, where it is for now assumed (and will be later verified) that $\partial_{xx}\widetilde{V}<0$. Namely,
\[
\pi^*_s(u;x) = \mathcal{N}\left(\frac{B \partial_x\widetilde{V}}{-D^2\partial_{xx}\widetilde{V}},\frac{\lambda_{SE}}{-D^2\partial_{xx}\widetilde{V}}\right).
\]
Therefore,
\[
\int_\mathbb{R} u \pi^*_s(u;x)du = \frac{B \partial_{x}\widetilde{V}}{-D^2\partial_{xx}\widetilde{V}} \,\,\text{ and } \int_\mathbb{R} u^2 \pi^*_s(u;x)du = \left(\frac{B \partial_{x}\widetilde{V}}{-D^2\partial_{xx}\widetilde{V}}\right)^2 + \frac{\lambda_{SE}}{-D^2\partial_{xx}\widetilde{V}}.
\]

\noindent Next, we introduce the ansatz 
\begin{eqnarray}\label{eqn:SE_ansatz}
\widetilde{V}(s,x) = -\frac{\eta_s}{2}(x-m_s)^2+\gamma_s,
\end{eqnarray}
for some $\eta_s>0$ and $\gamma_s$ to be appropriately defined.
 Then, $\partial_{x}\widetilde{V} = -\eta_s(x-m_s)$ and $\partial_{xx}\widetilde{V}=-\eta_s$, and thus,
\[
\int_\mathbb{R} u \pi^*_s(u;x)du = \frac{B(m_s-x)}{D^2}
\]
and 
\[
\int_\mathbb{R} u^2 \pi^*_s(u;x)du= \frac{B^2(x-m_s)^2}{D^4}+ \frac{\lambda_{SE}}{D^2\eta_s}.
\]
\noindent {\bf Step 2 (updating the mean-field information):} Denoting $\kappa = \frac{B}{D^2}$ and plugging in the forward equation for $p(s,x)$ yield
\begin{eqnarray}
\partial_s p(s,x) =&& -\partial_x \Big(\left((A+B \kappa)(m_s-x)\right)p(s,x)\Big) \nonumber\\
\vspace{10 mm}\nonumber
\\
&&\,\,\,\,+\frac{1}{2}D^2\partial_{xx}\left(\left(\kappa^2(x-m_s)^2+ \frac{\lambda_{SE}}{D^2\eta_s} \right)p(s,x)\right),\nonumber\\
\vspace{10 mm}\nonumber
\\
=&&\left(A+B \kappa\right)p(s,x)-\left((A+B \kappa)(m_s-x)\right) \partial_x p(s,x)\nonumber \\
\vspace{10 mm}\nonumber
\\
&&\,\,\,\,+\frac{1}{2}D^2\left(2\kappa^2p(s,x)+4\kappa^2(x-m_s)\partial_x p(s,x)\right)\nonumber\\
\vspace{10 mm}\nonumber
\\
&&\,\,\,\,+\frac{1}{2}D^2\left(\left(\kappa^2(x-m_s)^2+ \frac{\lambda_{SE}}{D^2\eta_s} \right)\partial_{xx}p(s,x)\right).\label{eq:inter_hjb}
\end{eqnarray}
\noindent {{\bf Step 3 (finding a fixed-point):}} Multiplying both sides of \eqref{eq:inter_hjb} by $x$ and integrating with respect to $x$ yields that
$d m_s = \left(\int x \partial_s p(s,dx)\right)ds = 0$.
Therefore, $m_s^*=  \mathbb{E}_{\xi \sim \nu}[\xi]=:m^*$. 

Furthermore, the HJB equation for $s\in[t,T)$ is  reduced to
\begin{eqnarray}\label{eq:hjb_reduce_se}
-\partial_s \widetilde{V}(s,x) =
&&\max_{\pi_s\in \PP(\mathbb{R})} \left(\left(A(m^*-x)+B \int_\mathbb{R} u \pi_s(u;x)du\right)\partial_{x}\widetilde{V}(s,x) \right.\\\vspace{10 mm}\nonumber
\\
&&\,\,\,\,\left.-\frac{Q}{2}(x-m^*)^2-\lambda_{SE} \int_\mathbb{R} \pi_s(u;x)\ln \pi_s(u;x)du+\frac{1}{2}D^2\left(\int_\mathbb{R} u^2 \pi_s(u;x)du\right)\partial_{xx}\widetilde{V}(s,x) \right), \nonumber
\end{eqnarray}
with $\widetilde{V}(T,x) = -\frac{\overline{Q}}{2}(m^*-x)^2$.
Plugging $\pi_s^*(u;x)$ and using ansatz \eqref{eqn:SE_ansatz} with $m_s=m^*$ into the above HJB give
\begin{eqnarray*}
\frac{\dot{\eta_s}}{2}(x-m^*)^2-\dot{\gamma}_s =&& -\left( A(m^*-x)+\frac{B^2(m^*-x)}{D^2}\right)\eta_s(x-m^*)\\
&-&\frac{Q}{2}(m^*-x)^2+\lambda_{SE} \ln\left(\sqrt{\frac{2\pi e \lambda_{SE}}{D^2\eta_s}}\right)\\
&-&\frac{1}{2}D^2\left(\frac{B^2(x-m^*)^2}{D^4}+ \frac{\lambda_{SE}}{D^2\eta_s}\right)\eta_s.
\end{eqnarray*}
Direct calculations imply
\begin{eqnarray}\label{eta33}
 \dot{\eta}_s =\left(2A+\frac{B^2}{D^2} \right)\eta_s-{Q},
 \end{eqnarray}
 with ${\eta}_T={\bar{Q}}$,  and 
 \begin{eqnarray}\label{gamma33}
\dot{\gamma}_s &=& - \frac{\lambda_{SE}}{2}\ln\left({\frac{2\pi  \lambda_{SE}}{D^2}}\right)+\frac{\lambda_{SE}}{2}\ln \eta_s
 \end{eqnarray}
 with  $\gamma_T=0$.  Then,
  \eqref{eta33} admits the unique solution 
  \begin{equation*}
  \eta_s = {\bar{Q}}\exp\left(-\left(2A+\frac{B^2}{D^2}\right)(T-s)\right)+{\frac{Q}{2A+\frac{B^2}{D^2}}\left(1- \exp\left({-\left(2A+\frac{B^2}{D^2}\right)(T-s)}\right)\right)},
  \end{equation*}
   from which it is easy to verify that $\eta_s>0$, since $A>0,Q>0$ and $\bar{Q}>0$, $s\in[t,T]$. Moreover, \eqref{gamma33} admits the unique solution
 $$\gamma_s =  \frac{\lambda_{SE}}{2}\ln\left({\frac{2\pi  \lambda_{SE}}{D^2}}\right)(T-s)-\int_s^T\frac{\lambda_{SE}}{2}\ln(\eta_z)dz.$$
 Consequently, one NE (optimal) policy takes the form 
\begin{equation*}
\pi_s^*(u;x) = \mathcal{N}\left(\frac{B(m^*-x)}{D^2},\frac{\lambda_{SE}}{D^2\eta_s}\right),
\end{equation*}
and the associated optimal state process is the unique solution of the SDE \eqref{eq:se_optimal_dynamics}.

\paragraph{Verification argument.}
\noindent The final step is to verify that $m^*$ is the mean state under policy \eqref{eqn:pi} and $V^*(t) := \mathbb{E}_{\xi\sim \nu}[\widetilde{V}(\xi,t)]=\mathbb{E}_{\xi\sim \nu}[-\frac{\eta_t}{2}(\xi-m^*)^2+\gamma_t]$ is the corresponding game value.

First, let us fix the mean field information as $m_s = m^*$, $s \in [t,T]$, and also fix the initial state $x\in\mathbb{R}$ and initial time $t \in [0,T]$. Let $\pi \in \mathcal{A}(x)$ and $X^{\pi}$ be the associated state process under $\pi$ solving
\[d X^{\pi}_s = \left(\int_\mathbb{R} (A(m^*-X^{\pi}_s)+Bu)\pi_s(u)du \right)ds + D\left(\sqrt{\int_\mathbb{R} u^2 \pi_s(u)du}\right) dW_s.
\]

Denote $\tilde{r}(x,\pi) = -\frac{Q}{2}(x-m^*)^2$, $\tilde{b}(x,\pi) = \int_{\mathbb{R}} (A(m^*-x)+Bu)\pi_s(u)du$, and $\tilde{\sigma}(x,\pi) = D\left(\sqrt{\int_\mathbb{R} u^2 \pi(u)du}\right)$. Further, define the stopping time $\tau_n^{\pi}:=\left\{s \ge t\,:\,\int_t^T\partial_{x}\widetilde{V}(t,X_{s}^{\pi})\tilde{\sigma}^2(X_{s}^{\pi},\pi_s)^2ds\ge n\right\}$, for $n \ge 1$. 
Then, It\^{o}'s formula yields
\begin{eqnarray*}
\widetilde{V}(T\wedge \tau_n^{\pi},X^{\pi}_{T \wedge \tau_n^{\pi}}) -\widetilde{V}(t,x) &=& \int_t^{T \wedge \tau_n^{\pi}}\left(\frac{1}{2}\partial_{xx}\widetilde{V}(s,X_s^{\pi}) \tilde{\sigma}^2(X_{s}^{\pi},\pi_s) +\partial_{x}\widetilde{V} (s,X_{s}^{\pi})\tilde{b} (X_{s}^{\pi},\pi_s)\right)ds \\
&&+ \int_t^{T \wedge \tau_n^{\pi}}\partial_{x}\widetilde{V}(s,X_{s}^{\pi})dW_s.
\end{eqnarray*}
Taking expectations,  using that $\widetilde{V}$ solves the HJB equation \eqref{eq:hjb_reduce_se}, and that $\pi$ is in general sub-optimal, we deduce that

\begin{eqnarray*}
&&\mathbb{E}\left[\widetilde{V}(T,X^{\pi}_{T \wedge \tau_n^{\pi}})\right] \\
\vspace{10 mm}\nonumber
\\
&=& \widetilde{V}(t,x) + \mathbb{E}\left [\int_t^{T \wedge \tau_n^{\pi}}\left(\frac{1}{2}\partial_{xx}\widetilde{V}(s,X_s^{\pi}) \tilde{\sigma}^2(X_{s}^{\pi},\pi_s) +\partial_{x}\widetilde{V} (s,X_{s}^{\pi})\tilde{b} (X_{s}^{\pi},\pi_s)\right)ds + \int_t^{T \wedge \tau_n^{\pi}}\partial_{x}\widetilde{V}(s,X_{s}^{\pi})dW_s\right]\\
&\leq& \widetilde{V}(t,x) - \mathbb{E}\left[\int_t^{T \wedge \tau_n^{\pi}} \left(\tilde{r}(X_{s}^{\pi},\pi_s)-\lambda \int_{\mathbb{R}} \pi_s(u)\ln \pi_s(u)du\right)ds\right].
\end{eqnarray*}
Standard calculations yield that $\mathbb{E}[\sup_{t\leq s \leq T}|X^{\pi}_s|^2]\leq N(1+x^2)e^{NT}$ for some constant $N>0$, which is independent of $n$. Sending $n \rightarrow \infty$ yields
\begin{eqnarray*}
\widetilde{V}(t,x) \geq  \mathbb{E}\left[\int_t^T \left(\tilde{r}(X_{s}^{\pi},\pi_s)-\lambda \int_{\mathbb{R}} \pi_s(u)\ln \pi_s(u)du\right)ds -\frac{\bar{Q}}{2}(X_T^{\pi}-m^*)^2\right],
\end{eqnarray*}
for each $x\in \mathbb{R}$ and $\pi \in \mathcal{A}$. Hence, $\widetilde{V}(t,x) \geq V^*(t,x)$, for all $x\in \mathbb{R}$.

On the other hand,  the right-hand of \eqref{hjb} is maximized for
\begin{equation}\label{eq:opt_pi_proof}
\pi_s^*(u;x) = \mathcal{N}\left(\frac{B(m^*-x)}{D^2},\frac{\lambda_{SE}}{D^2\eta_s}\right).
\end{equation}
Thus,
\[
\widetilde{V}(t,x) = \mathbb{E}\left[\int_0^T \left(\tilde{r}(X_{s}^{*},\pi_s)-\lambda \int_{\mathbb{R}} \pi^*_s(u;X_{s}^{*})\ln \pi^*_s(u;X_{s}^{*})du\right)ds -\frac{\bar{Q}}{2}(X_T^{*}-m^*)^2\right],\]
where $X_s^*$ is the controlled state process under policy \eqref{eq:opt_pi_proof}.

Next, let us show that for $s\in [t,T]$,
$$m^* = \mathbb{E}[X_s^*].$$ 

To this end, let
$K =-\left( A+\frac{B^2}{D^2}\right)$. Then,
$$d X_s^*=( KX_s^* -Km^*)ds + f(s,X_s^*,m^*)dW_s,$$
with  
$$f(s,x,m) = \left(\sqrt{\left(\frac{B}{D}(x-m)\right)^2+\frac{\lambda_{SE}}{\eta_s}}\right).$$
Therefore,
\begin{eqnarray*}
e^{- K(s-t)}X_s^* &=&\xi +\int_t^s e^{- K (z-t)}\left(-Km^*dz+ f(z,X_z^*,m^*)dW_z\right),
\end{eqnarray*}
and $e^{- K(s-t)}\mathbb{E}[X_s^*] = \mathbb{E}[\xi] +\Big(e^{- K(s-t)}-1\Big)m^*$. Hence,  $\mathbb{E}[X_s^*] =m^*$,  $s\in[t,T]$.

\end{proof}

\paragraph{NE Derivation of Game  {\bf (MFG-EE)}.} To ease the exposition, we drop the subscript EE.

\begin{proof}{Proof of Theorem \ref{thm:current_cross}}
For a given Markovian policy $\pi_s(u;x)$, the forward equation for $p(s,x)$, the density of $X_s$, $s \in [t,T]$ satisfies
\[
\partial_s p(s,x) =-\partial_x \left(\left(A({m}_s-x)+B \int_\mathbb{R} u \pi_s(u;x)du \right)p(s,x)\right) +\frac{1}{2}\partial_{xx}\left(p(s,x)\int_\mathbb{R} D^2u^2 \pi_s(u;x)du \right),
\]
with initial density $p(t,x) = \nu(x)$ and ${m}_s = \int x p(s,x)dx$.

\noindent {{\bf Step 1 (solving the control problem):}} Given fixed mean-field information $\{m_s\}_{s\in[0, T]}$, the HJB equation for the value function $\widetilde{V}(s,x,m)$ can be written as
\begin{eqnarray*}\label{hjb-6}
-\partial_s\widetilde{V}(s,x)=
&&\max_{\pi_s\in \mathcal{P}(\mathbb{R})} \left(\left(A(m_s-x)+B \int_\mathbb{R} u d\pi_s(u;x)\right) \partial_x\widetilde{V}(s,x) -\frac{Q}{2}(m_s-x)^2\right. \nonumber \\
\vspace{10 mm}\nonumber
\\
&&\,\,\,\,-\lambda_{CE} \int_\mathbb{R} \pi_s(u;x)\int \ln \alpha_s (u;x)\mu_s(dx)du-\lambda_{SE} \int_\mathbb{R} \pi_s(u;x)\ln \pi_s(u;x)du \nonumber  \\ 
\vspace{10 mm}\nonumber
\\
 &&\left. +\frac{1}{2}\left(\int_\mathbb{R} D^2u^2 \pi_s(u;x)du\right)\partial_{xx}\widetilde{V}(s,x)
\right), \nonumber
\end{eqnarray*}
with $\widetilde{V}(T,x) = -\frac{\bar{Q}}{2}(x-m_T)^2$. Recall that $\pi_s \in \mathcal{P}(U)$ if and only if \eqref{PU} holds.
The constrained maximization problem on the right hand side of \eqref{hjb} yields {\small
\begin{eqnarray*}
\pi_s^*(u;x) = \frac{\exp\left(\frac{1}{\lambda_{SE}}\left(-\frac{Q}{2}(x-m_s)^2+\frac{1}{2}D^2u^2\partial_{xx}\widetilde{V}+(A(m_s-x)+Bu)\partial_{x}\widetilde{V}-\lambda_{CE} \int \Big(\ln \alpha_s (u;x)\Big)\mu_s(dx)\right)\right)}{\int_\mathbb{R} \exp\left(\frac{1}{\lambda_{SE}}\left(-\frac{Q}{2}(x-m_s)^2+\frac{1}{2}D^2u^2\partial_{xx}\widetilde{V}+(A(m_s-x)+Bu)\partial_{x}\widetilde{V}-\lambda_{CE} \int \Big(\ln \alpha_s (u;x)\Big)\mu_s(dx)\right)\right)du}.
\end{eqnarray*}}
Next, we introduce the ansatz  for the population action distribution for the agent in state $y$,
\begin{eqnarray}\label{population_decision}
\alpha_s(u;y) = \mathcal{N}(u\,|\,H_s(y-m_s),L_s),
\end{eqnarray}
with some (to be defined) deterministic processes $H_s$ and $L_s>0$, $s\in[t,T]$. Then, $\alpha_s(u;y)$ is Gaussian with mean $H_s(y-m_s)$ and variance $L_s$. We stress that the Gaussian property of $\alpha_s(u;y)$ does not imply the Gaussian property of  the aggregated population action distribution $\widetilde{\alpha}(s) = \int \alpha_s(u;y)\mu_s(dy)$. 

In turn, $\ln \alpha_s(u;y) = -\frac{1}{2} \ln(2\pi L_s)-\frac{1}{2L_s}(u-H_s(y-m_s))^2$ and 
\begin{eqnarray*}
\int \mu_s(dy) \ln(\alpha_s(u;y)) &=& -\frac{1}{2}\ln(2\pi L_s) -\frac{1}{2L_s}\int(u-H_s(y-m_s))^2\mu_t(dy)\\
\vspace{10 mm}\nonumber
\\
 &=& -\frac{1}{2}\ln(2\pi L_s) -\frac{1}{2L_s}(u^2+H_s^2 \mbox{Var}(\mu_s)),
\end{eqnarray*}
with $\mbox{Var}(\mu_s) = \int x^2 \mu_s(dx) - (\int x \mu_s(dx))^2$.
Therefore, the optimal  policy  is \textit{Gaussian} with mean
$\frac{B \partial_{x}\widetilde{V}}{-D^2\partial_{xx}\widetilde{V}-\frac{\lambda_{CE}}{L_s}}$ and variance $\frac{\lambda_{SE}}{-D^2\partial_{xx}\widetilde{V}-\frac{\lambda_{CE}}{L_s}}$, namely,
\[
\pi_s^*(u;x) = \mathcal{N}\left(u\,\,\left\vert\,\,\frac{B \partial_{x}\widetilde{V}}{-D^2\partial_{xx}\widetilde{V}-\frac{\lambda_{CE}}{L_s}},\frac{\lambda_{SE}}{-D^2\partial_{xx}\widetilde{V}-\frac{\lambda_{CE}}{L_s}}\right.\right).\]

Let us  for now assume (and will verify later) that $-D^2\partial_{xx}\widetilde{V}-\frac{\lambda_{CE}}{L_s}>0$. In turn,
\[
\int_\mathbb{R} u \pi_s^*(u;x)du = \frac{B \partial_{x}\widetilde{V}}{-D^2\partial_{xx}\widetilde{V}-\frac{\lambda_{CE}}{L_s}} \] and 
\[ \int_\mathbb{R} u^2 \pi_s^*(u;x)du = \left(\frac{B \partial_{x}\widetilde{V}}{-D^2\partial_{xx}\widetilde{V}-\frac{\lambda_{CE}}{L_s}}\right)^2 + \frac{\lambda_{SE}}{-D^2\partial_{xx}\widetilde{V}-\frac{\lambda_{CE}}{L_s}}.
\]
Next, consider the ansatz 
\begin{eqnarray}\label{eq:EE_ansatz}
\widetilde{V}(s,x) =- \frac{\eta_s}{2}(x-m_s)^2+\gamma_s.
\end{eqnarray}
In turn, $\partial_{x}\widetilde{V} = -\eta_s(x-m_s)$ and $\partial_{xx}\widetilde{V}=-\eta_s$, together with
\[
\int_\mathbb{R} u \pi_s^*(u;x)du = \frac{B\eta_s(m_s-x)}{D^2\eta_s-\frac{\lambda_{CE}}{L_s}},
\]
and 
\[
\int_\mathbb{R} u^2 \pi_s^*(u;x)du  = \left(\frac{B\eta_s(m_s-x)}{D^2\eta_s-\frac{\lambda_{CE}}{L_s}}\right)^2+ \frac{\lambda_{SE}}{D^2\eta_s-\frac{\lambda_{CE}}{L_s}}.
\]
\noindent {{\bf Step 2 (updating the mean-field information):}} Denoting $\kappa_s := \frac{B\eta_s}{D^2\eta_s-\frac{\lambda_{CE}}{L_s}}$, $s\in[t,T]$, and plugging in the forward equation for $p(s,x)$, we deduce that
\begin{eqnarray}
\partial_s p(s,x) =&& -\partial_x \left((A+B \kappa_s)(m_s-x)p(s,x)\right)\nonumber \\
\vspace{10 mm}\nonumber
\\
&&\,\,\,\,+\frac{1}{2}D^2\partial_{xx}\left(\left(\kappa_t^2(x-m_s)^2+ \frac{\lambda_{SE}}{D^2\eta_s-\frac{\lambda_{CE}}{L_s}} \right)p(s,x)\right)\nonumber\\
\vspace{10 mm}\nonumber
\\
=&&\left(A+B \kappa_s\right)p(s,x)-\left((A+B \kappa_s)(m_s-x)\right)\partial_x p(s,x) \nonumber\\
\vspace{10 mm}\nonumber
\\
&&\,\,\,\,+\frac{1}{2}D^2\left(2\kappa_s^2p(s,x)+4\kappa^2_s(x-m_s)\partial_x p(s,x)\right)\nonumber\\
\vspace{10 mm}\nonumber
\\
&&\,\,\,\,+\frac{1}{2}D^2\left(\left(\kappa_s^2(x-m_s)^2+ \frac{\lambda_{SE}}{D^2\eta_s-\frac{\lambda_{CE}}{L_s}} \right)\partial_{xx} p(s,x)\right).\label{hjb_tmp}
\end{eqnarray}
Multiplying both sides of \eqref{hjb_tmp} by $x$ and integrating with respect to $x$ give $\partial_s m_s = \int x \partial_s p(s,dx) = 0$ and, thus, $m^*_s = m^* = \mathbb{E}_{\xi\sim\nu}[\xi]$, for $ s \in [t  , T]$. Hence, the HJB equation reduces to
\begin{eqnarray*}
-\partial_s \widetilde{V}(s,x)=
&\max_{\pi_s\in \PP(\mathbb{R})} \left(\left(A(m^*-x)+B \int_\mathbb{R} u d\pi_s(u)\right)\partial_x\widetilde{V}(s,x) -\frac{Q}{2}(m^*-x)^2\right.\\
\vspace{10 mm}\nonumber
\\
&\,\,\,\,\left.-\lambda_{SE} \int_\mathbb{R} \pi_s(u)\ln \pi_s(u)du-{\lambda_{CE} \int_\mathbb{R} \pi_s(u)\int \mu_s (dy)\ln {\alpha}_s(u;y)dydu}\right.\\
\vspace{10 mm}\nonumber
\\
&\left.+\frac{1}{2}D^2\int_\mathbb{R} u^2 \pi_s(u)du\partial_{xx}\widetilde{V}(s,x)\right).
\end{eqnarray*}
Plugging $\pi_s^*(u;x)$ and using ansatz \eqref{eq:EE_ansatz} with $m_s=m^*$ for the above HJB, we obtain

\begin{eqnarray*}
\frac{\dot{\eta_s}}{2}(x-m^*)^2-\dot{\gamma}_s =&-& \left( A(m^*-x)+B\frac{B\eta_s(m^*-x)}{D^2\eta_s-\frac{\lambda_{CE}}{L_s}}\right)\eta_s(x-m^*)+\lambda_{SE} \ln\left(\sqrt{\frac{2\pi e \lambda_{SE}}{D^2\eta_s-\frac{\lambda_{CE}}{L_s}}}\right)\\
&-&\frac{1}{2}D^2\left(\left(\frac{(B\eta_s)(x-m^*)}{D^2\eta_s-\frac{\lambda_{CE}}{L_s}}\right)^2+ \frac{\lambda_{SE}}{D^2\eta_s-\frac{\lambda_{CE}}{L_s}}\right)\eta_s-\frac{Q}{2}(x-m^*)^2\\
&+& \frac{\lambda_{CE}}{2}\ln(2\pi L_s) +\frac{\lambda_{CE}}{2L_s}\left(\left(\frac{B\eta_s(x-m^*)}{D^2\eta_s-\frac{\lambda_{CE}}{L_s}}\right)^2+\frac{\lambda_{SE}}{D^2\eta_s-\frac{\lambda_{CE}}{L_s}}\right)\\
&+& \frac{\lambda_{CE}}{2L_s}H_s^2 {\mbox{Var}}(\mu_s).
\end{eqnarray*}
Direct calculations yield 
\begin{eqnarray}
\dot{\eta}_s &=& 2A\eta_s+\frac{\left(B\eta_s\right)^2}{D^2\eta_s-\frac{\lambda_{CE}}{L_s}}-Q,\label{eta-4}\\
\vspace{10 mm}\nonumber
\\
\dot{\gamma}_s &=&\frac{\lambda_{SE}}{2}-\lambda_{SE} \ln \left(\sqrt{\frac{2\pi e \lambda_{SE}}{D^2\eta_s-\frac{\lambda_{CE}}{L_s}}} \right)-\frac{\lambda_{CE}}{2} \ln \left(2\pi L_s\right)-\frac{\lambda_{CE}}{2L_s}H_s^2 {\mbox{Var}}(\mu_s).\label{gamma-4}
\end{eqnarray}
\noindent {{\bf Step 3 (finding a fixed-point):}} Setting
\[
L_s = \frac{\lambda_{SE}}{D^2\eta_s-\frac{\lambda_{CE}}{L_s}}\,\,\, \textit{ and }\,\,\,
H_s = \frac{-B \eta_s}{D^2\eta_s-\frac{\lambda_{CE}}{L_s}},
\]
we deduce that
\[
H_s = -\frac{B}{D^2}\frac{\lambda_{SE}+\lambda_{CE}}{\lambda_{SE}}\,\,\, \textit{ and }\,\,\,
L_s = \frac{\lambda_{SE}+\lambda_{CE}}{D^2\eta_s},
\]
with
\begin{eqnarray}
\dot{\eta}_s &=& 2A\eta_s +\frac{B^2\eta_s}{D^2}\frac{\lambda_{SE}+\lambda_{CE}}{\lambda_{SE}}-Q,\label{eta3}\\
\vspace{10 mm}\nonumber
\\
\dot{\gamma}_s &=& -\frac{\lambda_{SE}+\lambda_{CE}}{2} \ln \left( \frac{2\pi(\lambda_{SE}+\lambda_{CE})}{D^2\eta_s}\right)-\frac{\lambda_{CE}}{2}\frac{B^2\eta_s(\lambda_{SE}+\lambda_{CE})}{D^2\lambda_{SE}^2}\mbox{Var}(\mu_s).\label{eq:(b)}
\end{eqnarray}
Consequently,
\[
\pi_s^*(u;x) = \mathcal{N}\left(u\,\,\left|\,\,\frac{B}{D^2}{\frac{\lambda_{SE}+\lambda_{CE}}{\lambda_{SE}}}(m^*-x), \frac{\lambda_{SE}+\lambda_{CE}}{D^2\eta_s}\right.\right).
\]
Denote $-K =\left( A+B\frac{B}{D^2}\right){\frac{\lambda_{SE}+\lambda_{CE}}{\lambda_{SE}}}$ and 
$$f(s,x,m) = \left(\sqrt{\left(\frac{B}{D}{\frac{\lambda_{SE}+\lambda_{CE}}{\lambda_{SE}}}(x-m)\right)^2+\frac{\lambda_{SE}+\lambda_{CE}}{\eta_s}}\right).$$
In turn,
$$d X_s^*= K(X_s^*-m^*)ds + f(s,X_s^*,m^*)dW_s,$$
and 
\begin{eqnarray*}
d(e^{-K(s-t)}X_s^*) &=& -Ke^{-K(s-t)}X_s^*ds+e^{-K(s-t)}dX_s^* \\
\vspace{5 mm}\nonumber
\\
&=&-Ke^{-K(s-t)}X_s^*ds+e^{-K(s-t)}(( KX_s^* -Km^*)ds+ f(s,X_s^*,m^*)dW_s)\\
\vspace{5 mm}\nonumber
\\
&=&e^{-K(s-t)}(-Km^*ds+ f(s,X_s^*,m^*)dW_s).
\end{eqnarray*}
Therefore,
\begin{eqnarray}\label{eq:controlled_dynamics}
e^{- K(s-t)}X_s^* &=&\xi +\int_t^s e^{- K (z-t)}(-Km^*dz+ f(z,X_z^*,m^*)dW_z),
\end{eqnarray}
and
\[
e^{- 2K(s-t)}\mbox{Var}[X^*_s] = \mbox{Var}[\xi]+\E\left[\left(\int_t^{s} e^{- K (z-t)}f(z,X^*_z,m^*)dW_z\right)^2\right].
\]
By It\^o's isometry, 
\begin{eqnarray*}
&&\E\left[\left(\int_t^{s} e^{- K (z-t) }f(z,X^*_z,m^*)dW_z\right)^2\right] = \E\left[\int_t^{s} e^{-2 K (z-t) }f^2(z,X^*_z,m^*)dz\right]\\
&=&\E\left[\int_t^{s} e^{- 2K (z-t) }\left(\left(\frac{B}{D}{\frac{\lambda_{SE}+\lambda_{CE}}{\lambda_{SE}}}(X_z^*-m^*)\right)^2+\frac{\lambda_{SE}+\lambda_{CE}}{\eta_z}\right)dz\right]\\
&=&\int_t^{s} e^{- 2K (z-t) }\left(\left(\frac{B}{D}{\frac{\lambda_{SE}+\lambda_{CE}}{\lambda_{SE}}}\right)^2\mbox{Var}[X_z^*]+\frac{\lambda_{SE}+\lambda_{CE}}{\eta_z}\right)dz.
\end{eqnarray*}

{
\noindent Let $$y(s) = e^{- 2K(s-t)}\mbox{Var}[X^*_s],$$  $$M = \left(\frac{B}{D}{\frac{\lambda_{SE}+\lambda_{CE}}{\lambda_{SE}}}\right)^2 \,\,{\rm and}\quad b(s) = e^{- 2K (s-t) }\frac{\lambda_{SE}+\lambda_{CE}}{\eta_s}.$$ 
Thus,
\[
y(s) = e^{M(s-t)} \left( y(t)+\int_t^s e^{-M(z-t)} b(z) dz \right),
\]
and
\[
e^{-2K(s-t)} \mbox{Var} (X_s^*) = e^{M(s-t)} \left(\mbox{Var}(\xi) +\int_t^s e^{-(M+2K)(z-t)}\frac{\lambda_{SE}+\lambda_{CE}}{\eta_z}dz \right).
\]
Therefore,
\begin{eqnarray}\label{variance}
\mbox{Var}[X^*_s] &=& e^{ (2K+M)(s-t)} \mbox{Var}[\xi] + e^{ 2K(s-t)} \int_t^s e^{M (s-z)} b(z) dz\nonumber\\
&=&e^{ (2K+M)(s-t)} \mbox{Var}[\xi] + \int_t^s e^{(M+2K) (s-z)}\frac{\lambda_{SE}+\lambda_{CE}}{\eta_z} dz.
\end{eqnarray}}
Assume for the moment that $\eta_s>0$ for $s\in[t,T]$. Then, $\mbox{Var}[X^*_s]$ is well-defined.
%
Hence \eqref{eq:(b)} reduces to
 \begin{eqnarray}\label{gamma3}
\dot{\gamma}_s &=& - \frac{\lambda_{SE}+\lambda_{CE}}{2}\ln\, {\frac{2\pi  (\lambda_{SE}+\lambda_{CE})}{D^2}} \nonumber\\
&& +\frac{\lambda_{SE}+\lambda_{CE}}{2}\ln \eta_s-\frac{\lambda_{CE}}{2}\frac{B^2\eta_s(\lambda_{SE}+\lambda_{CE})}{D^2\lambda_{SE}^2}\kappa_s,
 \end{eqnarray}
 with $\gamma_T=0$, where $\kappa_s$, $s\in[t,T]$,
 \begin{eqnarray}
\kappa_s:=e^{ (M+2K)(s-t)} \mbox{Var}[\xi] + \int_t^s e^{(M+2K) (s-z)}\frac{\lambda_{SE}+\lambda_{CE}}{\eta_z} dz
 \end{eqnarray}
with 
\begin{eqnarray*}
 K =-\left( A+B\frac{B}{D^2}\right){\frac{\lambda_{SE}+\lambda_{CE}}{\lambda_{SE}}} \,\, \text{ and }\,\,
M = D^2\left(\frac{B}{D^2}{\frac{\lambda_{SE}+\lambda_{CE}}{\lambda_{SE}}}\right)^2.
\end{eqnarray*}

\noindent Therefore, equation \eqref{eta3} admits the unique solution 
\begin{eqnarray*}
\eta_s &=& {\bar{Q}}\exp\left(-\left(2A+\frac{B^2}{D^2}\frac{\lambda_{SE}+\lambda_{CE}}{\lambda_{SE}}\right)(T-s)\right)\\
&&+\frac{Q}{2A+\frac{B^2}{D^2}\frac{\lambda_{SE}+\lambda_{CE}}{\lambda_{SE}}}\left(1-\exp\left(-\left(2A+\frac{B^2}{D^2}\frac{\lambda_{SE}+\lambda_{CE}}{\lambda_{SE}}\right)(T-s)\right)\right).
\end{eqnarray*}
We easily deduce that $\eta_s>0$, $s\in[t,T]$, since $\bar{Q}>0, Q>0$ and $A>0$.

Moreover, equation \eqref{gamma3} admits the (unique) solution
\begin{eqnarray*}
\gamma_s &=&  \frac{\lambda_{SE}+\lambda_{CE}}{2}(T-s)\ln\,{\frac{2\pi ( \lambda_{SE}+\lambda_{CE})}{D^2}}\,-\frac{\lambda_{SE}+\lambda_{CE}}{2}\int_s^T\,\ln\,\eta_z\,dz\\
&& +\int_s^T\frac{\lambda_{CE}}{2}\frac{B^2\eta_z(\lambda_{SE}+\lambda_{CE})}{D^2\lambda_{SE}^2}\kappa_z dz.
\end{eqnarray*}
We then obtain that the associated optimal feedback policy is given by
\begin{equation}
\pi_s^*(u;x) = \mathcal{N}\left(\frac{\lambda_{SE}+\lambda_{CE}}{\lambda_{SE}}\frac{B(m^*-x)}{D^2},\frac{\lambda_{SE}+\lambda_{CE}}{D^2\eta_s}\right),
\end{equation}
and the optimal controlled state process is the unique solution of the SDE \eqref{eq:ee_optimal_dynamics}.
The verification is similar to the verification of Theorem \ref{thm:current_shannon} and is therefore omitted.
\end{proof}

\paragraph{Proof of Corollary \ref{cor:random_policy}}
\begin{proof}{}
Let $\hatK := -\left(A+B\hatM\right)$ and $f\left(m_s,X_s,\hatM,\hatV_s\right) := D\sqrt{\hatM^2(m_s-X_s)^2 +\hatV_s}$, $s\in[t,T]$. Then, under policy $\widehat{\pi}$ given in \eqref{eqn:hatcontrol},
\begin{eqnarray*}
d X^{\widehat{\pi}}_s &=& \left(\left(A+B\hatM\right)m_s - \left(A+B\hatM \right)X^{\widehat{\pi}}_s \right) ds + f\left(m_s,X^{\widehat{\pi}}_s,\hatM,\hatV_s\right) dW_s\\
&=& -\widehat{K} m_sds +KX^{\widehat{\pi}}_s ds + f  dW_s.
\end{eqnarray*}

\noindent Using that
$d(e^{- \hatK s}X^{\widehat{\pi}}_s)=e^{-\hatK s}\left(-\hatK m_sds + f\,dW_s\right)$,
we have
\begin{eqnarray}\label{eqn:intermediate}
e^{- \hatK (s-t)}X^{\widehat{\pi}}_s &=&\xi +\int_t^s e^{- \hatK (z-t)}(-\hatK m_zdz+ f\,dW_z).
\end{eqnarray}
Hence,
\begin{eqnarray*}
\E [X^{\widehat{\pi}}_s] &=&e^{\hatK (s-t)}\E [\xi] - \int_t^s e^{\hatK (s-z)}\hatK m_z\,dz.
\end{eqnarray*}
Let $\widehat{m}_s :=\E [X^{\widehat{\pi}}_s]$. From \eqref{eqn:intermediate} and routine calculations, we deduce that
\begin{eqnarray*}
e^{- 2\hatK (s-t)}\left(X^{\widehat{\pi}}_s\right)^2
&=&\xi^2 +2\xi\,\int_t^s e^{- \hatK (z-t)}(-\hatK m_zdz+ f\,dW_z)+\left(\int_t^s e^{- \hatK (z-t)}\hatK m_zdz\right)^2\\
&&+\left(\int_t^s e^{- \hatK (z-t)} f\,dW_z\right)^2 -2 \int_t^s e^{- \hatK (z-t)}\hatK m_zdz\,\int_t^s e^{- \hatK (z-t)} f\,dW_z.
\end{eqnarray*}
Hence,
\begin{eqnarray}\label{eqn:square_expecation}
&&e^{- 2\hatK (s-t)}\E \left[\left(X^{\widehat{\pi}}_s\right)^2\right] \nonumber\\
&&=\E[\xi^2] -2\E[\xi]\,\int_t^s e^{- \hatK (z-t)}\hatK m_zdz+\left(\int_t^s e^{- \hatK (z-t)}\hatK m_zdz\right)^2+\E\left(\int_t^s e^{- \hatK (z-t)} f\,dW_z\right)^2.
\end{eqnarray}
By It\^{o}'s isometry, 
\begin{eqnarray}\label{eqn:itoisometry}
&&\E\left(\int_t^s e^{- \hatK (z-t)} f\,dW_z\right)^2\nonumber 
=\E\left[ \int_t^s  e^{- 2\hatK (z-t)} f^2\,d z\right]\nonumber\\
 &=& \E\left[ \int_t^s e^{-2 \hatK (z-t)}D^2\left(\hatM^2\left(X^{\widehat{\pi}}_z-m_z\right)^2 +\hatV_z\right)d z\right]\nonumber\\
  &=& \E\left[ \int_t^s e^{-2 \hatK (z-t)}D^2\hatM^2\left(m_z^2-2m_zX^{\widehat{\pi}}_z+\left(X^{\widehat{\pi}}_z\right)^2\right)d z +\int_t^se^{-2 K (z-t)}D^2\hatV_zd z\right]\nonumber\\
  &=& \E\left[ \int_t^s e^{- 2\hatK (z-t)}D^2\hatM^2\left(X^{\widehat{\pi}}_z\right)^2d z \right]+\int_t^s e^{- 2\hatK (z-t)}D^2\left(\hatM^2m_z^2-2\hatM^2m_z \widehat{m}_z +\hatV_z \right)d z\nonumber\\
  &=&  D^2\hatM^2\int_t^s e^{- 2 \hatK (z-t)}\E\left[\left(X^{\widehat{\pi}}_z\right)^2\right]d z +\int_t^se^{- 2\hatK (z-t)}D^2\left(\hatM^2m_z^2-2\hatM^2m_z \widehat{m}_z +\hatV_z \right)d z
\end{eqnarray}
Combining \eqref{eqn:square_expecation} and \eqref{eqn:itoisometry} yields
\begin{eqnarray*}
e^{- 2K(s-t)}\E \left[\left(X^{\widehat{\pi}}_s\right)^2\right] &=&\E[\xi^2] -2\E[\xi]\,\int_t^s e^{- \hatK (z-t)}\widehat{K}m_z dz+\left(\int_t^s e^{- \hatK (z-t)}\hatK m_zdz\right)^2\\
&&  + D^2\hatM^2\int_t^s e^{-2 \hatK (z-t)}\E \left[\left(X^{\widehat{\pi}}_z\right)^2\right]d z +\int_t^se^{- 2\hatK (z-t)}D^2\left(\hatM^2m_z^2-2\hatM^2m_z \widehat{m}_z +\hatV_z \right)d z.
\end{eqnarray*}
Letting $y_s := e^{- 2\hatK (s-t)}\E \left[\left(X^{\widehat{\pi}}_s\right)^2\right]$ for $s\in[t,T]$, we have
\[
y_s -y_t= b_s +  D^2\hatM^2\int_t^s y_z dz, 
\]
where, 
\[
b_s := -2\E[\xi]\,\int_t^s e^{- \hatK (z-t)}\hatK m_zdz+\left(\int_t^s e^{- \hatK (z-t)}\hatK m_zdz\right)^2 +\int_t^se^{-2 \hatK (z-t)}D^2\left(\hatM^2m_z^2-2\hatM^2m_z \widehat{m}_z +\hatV_z \right)d z,
\]
with $b_t=0$. {Therefore,
\[
\int_t^s \dot{y}_z dz  = \int_t^s \dot{b}_z dz +  D^2\hatM^2\int_t^s y_z dz
\]
and, thus,
\[
y_s =  e^{ D^2\hatM^2(s-t)}\left(y_t+ \int_t^s e^{- D^2\hatM^2(z-t)}\dot{b}(z)dz\right).\]
Finally, for $s\in[t,T]$,
\[
\E\left[\left(X^{\widehat{\pi}}_s\right)^2\right] = e^{\left(2\hatK +D^2\hatM^2\right)(s-t)}\left(\E[\xi^2] + \int_t^s e^{- D^2\hatM^2(z-t)}\dot{b}(z)dz\right)
\]
with
\begin{eqnarray*}
\dot{b}(s) &=&  -2\E[\xi]\, e^{- \hatK (s-t)}\hatK m_s+\left(\int_t^s e^{- \hatK (z-t)}\hatK m_zdz\right) e^{- \hatK (s-t)}\hatK m_s \\
&&+e^{-2 \hatK (s-t)}D^2\left(\hatM^2m_s^2-2\hatM^2m_s \widehat{m}_s +\hatV_s \right).
\end{eqnarray*}
}
The rest of the proof follows easily.
\end{proof}

\section{Experiment}\label{sec:experiment}
We now demonstrate how  the theoretical results of Theorems \ref{thm:current_shannon} and \ref{thm:current_cross} can be used to design algorithms for  MFG with learning.  The experiment aims to highlight
\begin{itemize}
    \item how  entropy regularization  helps  to ``explore optimally'' in a game with learning, and especially in improving the speed of convergence to the NE, and 
    \item how the agent manages to  eventually learn the optimal scheduling of the exploration and, in particular, the time-dependent variances  (as  in \eqref{eqn:pi} and \eqref{eqn:pi2}) over a finite time horizon. 
\end{itemize}
Throughout this section, the experiment  is  with the inclusion of Shannon entropy only, as
 the case with the additional cross-entropy may be studied in a  similar fashion.
\subsection{Set-up}
The algorithm design is with discrete time steps $s=0,1,2,\cdots,N,$ where  $\delta = \frac{T}{N}$ is the step-size. 
According to Theorem \ref{thm:current_shannon} and Corollary \ref{cor:random_policy}, it suffices to focus on a considerably smaller class of policies of form 
\[
\widehat{\pi}_s \sim \mathcal{N}\left(\hatM(x_s-m_s),\hatV_s \right),
\]
which can be fully characterized by the mean state process $m:=\{m_s\}_{s=0}^N$ and $\hatR:=(\hatM,\hatV)$, with $\hatV:=\{\hatV_s\}_{s=0}^N$. 
We, then, consider the discrete-time LQ-MFG problem 
\begin{equation}\label{objective:discrete-time}
    J \left(\hatR,m\right):= \mathbb{E}\left[\sum_{s=0}^{N-1}\left(-\frac{Q}{2}(x_s-m_s)^2+\lambda_{SE}\mathcal{H}_{SE}(\pi_s)\right)\delta -\frac{\bar{Q}}{2}(X_N-m_N)^2\right],
\end{equation}
where, for $s=0,1,\cdots,N-1$,
\begin{equation}\label{dynamics:discrete-time}
    x_{s+1} = x_s + \left(\int_\mathbb{R} (A(m_s-X_s)+B u)\pi_s(u)du\right)\delta +\left(D\sqrt{\int_\mathbb{R} u^2\pi_s(u)du} \right)\Delta W_s,\ x_0=\xi\sim\nu.
\end{equation}
Here,  $\Delta W_s$ are i.i.d $\mathcal{N}(0,\delta)$ random variables and  $\nu$ is the distribution of the initial state $\xi$.

\subsection{Mean Field Policy Gradient with Exploration}
Recall that in the learning setting, the model parameters $A$, $B$, $D$, $Q$, and  $\bar{Q}$ are assumed to be unknown to the agent.
She only has  access to the   {\it simulated} reward function 
\[
\widehat{j} \left(\widehat{R},m\right) :=\sum_{s=0}^{N-1}\Big((x_s-m_s)^2+\lambda_{SE}\mathcal{H}_{SE}(\pi_s)\Big)\delta -\frac{Q}{2}(X_N-m_N)^2,
\]
which is associated with a {\it single} trajectory $\{x_s\}_{s=0}^N$ under the policy characterized by $\widehat{R}$ and  the mean state process $m=\{m_s\}_{s=0}^N$. Note, however, that this assumption is weaker than being able to  observe $ J \left(\hatR,m\right)$ defined in \eqref{objective:discrete-time}, as $ J \left(\hatR,m\right)$ involves calculating an expectation and requiring observing infinite number of samples.

\begin{algorithm}[ht!]
\caption{\textbf{Mean Field Policy Gradient with Exploration}}
\label{alg:MFP_mf}
\begin{algorithmic}[1]
    \STATE \textbf{Input}: Initial beliefs $m^0 := \{m_s^0\}_{s=0}^{N}$, distribution of initial policy $ \mathcal{D}$, number of trajectories $n$, smoothing parameter $r$,  learning rate $\eta$.
       \FOR {$k \in \{1,\ldots,K\}$}
       \STATE Sample initial policy $\hatR^0 \sim \mathcal{D}$.
        \FOR {$i\in\{0, \ldots, I\}$}
             \FOR {$j\in\{1, \ldots, n\}$}
           \STATE Sample policy $\hatR^{i,j} =\hatR^{i} +U^{i,j}$ where $U^{i,j}\in \mathbb{R}^{N+2}$ is drawn uniformly at random over matrices such that $\|U^{i,j}\|_F=r$, where $\|\cdot\|_F$ is the Frobenius norm.
           \STATE Denote  {$\widehat{j}\left(\hatR^{i,j},m^{k-1}\right)$} as the {single trajectory cost} with policy $\hatR^{i,j}$ starting from $x^{i,j}_0 \sim \mathcal{\nu}$ under fixed mean state $m^{k-1}$.
        \ENDFOR
    \STATE Obtain  the estimate of $\nabla J(\hatR^i,m^{k-1})$:    \begin{eqnarray}\label{eqn:biased_estimate_gradient}
\widehat{\nabla J(\hatR^i,m^{k-1})} = \frac{1}{n}\sum_{j=1}^n \frac{1}{r^2}\,{\widehat{j}\left(\hatR^{i,j},m^{k-1}\right)}\, U^{i,j}.
 \end{eqnarray}
           
           \STATE Perform policy gradient descent step:
           \begin{eqnarray}\label{eqn:model_free_policy_update_2}
\hatR^{i+1} = \hatR^{i} -\eta \widehat{\nabla J(\hatR^i,m^{k-1})}.
\end{eqnarray}
          \ENDFOR
          \STATE Update mean field information $m^{k}=\{m_s^k\}_{s=0}^N$ assuming all agents follow policy $\hatR^{I+1}$.
        \ENDFOR
\end{algorithmic}
\end{algorithm}

The algorithm has the following key elements:
\begin{itemize}
    \item  {\it Information adaptiveness}. In each outer iteration $k$, $k=1,2,\cdots,K$, the representative agent can improve her decision (lines 4-11) based on the mean field information $m^{k-1}$ from the previous outer iteration $k-1$. This implies that she has access to a simulator with which she may exercise different policies when other agents keep applying the same policy from the previous outer iteration. This is a standard assumption, see, for example,  \cite{FYCW2019,GHXZ2019}. Once the agent stops improving her policy, the mean field information is updated assuming that all agents follow the same improved policy (line 12).  In the RL literature, this procedure is sometimes called  {\it
  fictitious play} \cite{cardaliaguet2017learning}.
    \item {\it Agent update.} Within each outer iteration $k$ under a fixed mean field information $m^{k-1}$, the agent will update her estimation of the optimal policy $\widehat{R}$ for $I$ rounds (lines 4-11). Each round corresponds to one gradient descent step (line 10) and requires $n$ samples of the simulated reward function (line 7) associated with the perturbed  version of $\widehat{R}^{i}$ (line 6). 
    
    \item The gradient term $\nabla J(\hatR^i,m^{k-1})$ in \eqref{eqn:biased_estimate_gradient} is estimated  using a {\it zeroth-order optimization approach} (line 9).
    That is, the agent only has query access to a sample of the reward function $\widehat{j}(\cdot)$ at input points $(R,m)$, without querying the gradients and higher order derivatives of $\widehat{j}(\cdot)$. 
  Moreover, to avoid the issue of ill-definedness of $\mathbb{E}_{U \sim \mathcal{N} (0,\sigma^2I)}[J(\hatR+U,m)]$ with a Gaussian smoothing, we choose $\mathbb{S}_r$ by smoothing over the sphere of a ball; hence, step \eqref{eqn:biased_estimate_gradient} in Algorithm \ref{alg:MFP_mf} is to find, for a given $m$, a bounded and biased estimate $\widehat{\nabla J(\hatR,m)} $ of $\nabla J(\hatR,m) $.
\end{itemize}

\subsection{Results}
\paragraph{Model set-up.} We take $T=0.1$, $\delta = 0.02$ {(hence $N=\frac{T}{\delta}=5$)}, $A = 2.0$, $B = 3.0$, $D = 2.0$, $Q = 3.0$, $\bar{Q} = 2.0$, $\mathbb{E}[\xi]=0.1$, and $\mathbb{E}[\xi^2]=1$.
\paragraph{Experiment set-up.} We set $r=0.01$ and $\eta = 0.05$, $m_s^0=0.0$ ($s=0,1,\cdots, N$), $\widehat{M}^0 \sim \mathcal{N}(0.5,1)$, and $\widehat{\sigma}_s^0 \sim \mathcal{N}(0.5,0.1)$ ($s=0,1,\cdots, N-1$), with $K=10$, $n=50$, and $I=400$.
\paragraph{Performance evaluation.}
Given policy $\hatR$ and mean field information $m$,  define the {\it relative error} between $(\hatR,m)$ and the mean field solution  $(R^*,m^*)$ of problem \eqref{objective:discrete-time}-\eqref{dynamics:discrete-time} as
\begin{eqnarray}
\rm{Err}(\hatR,m) := \frac{|J(\hatR,m)-J({R}^*,m^*)|}{|J({R}^*,m^*)|}.
\end{eqnarray}

\begin{figure}[H]
\centering
\begin{subfigure}[t]{.24\textwidth} 
\centering
\includegraphics[width=1.0\linewidth]{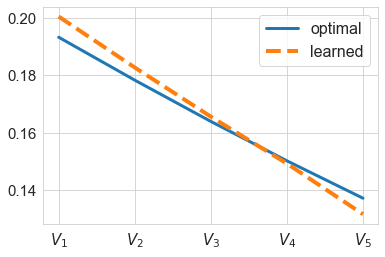}
\caption{\label{fig:lambda1_V}$\{V_s\}_{s=0}^N$}
\end{subfigure}
\begin{subfigure}[t]{.24\textwidth} 
\centering
\includegraphics[width=1.0\linewidth]{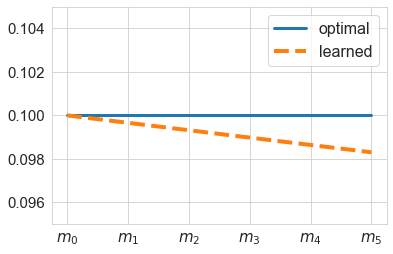}
\caption{\label{fig:lambda1_m}$\{m_s\}_{s=0}^N$}
\end{subfigure}
\begin{subfigure}[t]{.24\textwidth} 
\centering
\includegraphics[width=1.0\linewidth]{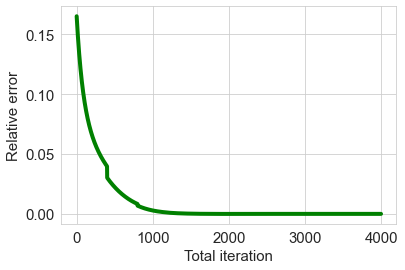}
\caption{Relative error in total iterations.}
\end{subfigure}
\begin{subfigure}[t]{.24\textwidth} 
\centering
\includegraphics[width=1.0\linewidth]{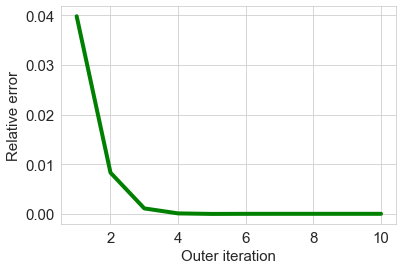}
\caption{\label{fig:lambda1_outer}Relative error in outer iterations}
\end{subfigure}
\caption{Performance of the algorithm when $\lambda_{SE}=1.0$. (True $M=0.75$ and learned $\widehat{M} = 0.732$.)}
\end{figure}

	\begin{figure}[H]
\centering
\begin{subfigure}[t]{.24\textwidth} 
\centering
\includegraphics[width=1.0\linewidth]{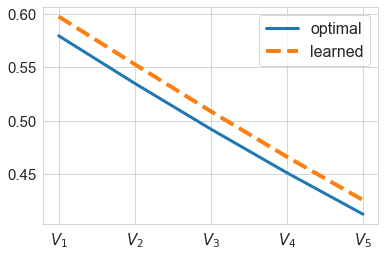}
\caption{\label{fig:lambda3_V}$\{V_s\}_{s=0}^N$}
\end{subfigure}
\begin{subfigure}[t]{.24\textwidth} 
\centering
\includegraphics[width=1.0\linewidth]{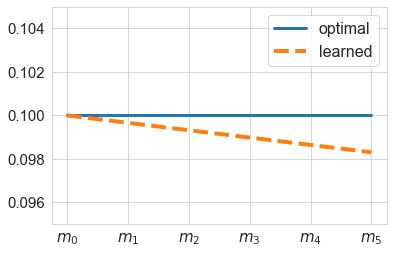}
\caption{\label{fig:lambda3_m}$\{m_s\}_{s=0}^N$}
\end{subfigure}
\begin{subfigure}[t]{.24\textwidth} 
\centering
\includegraphics[width=1.0\linewidth]{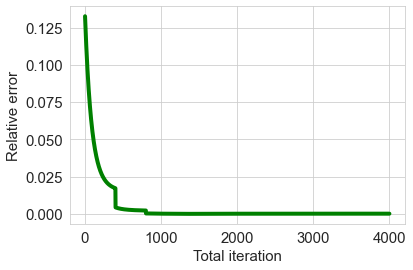}
\caption{Relative error in total iterations.}
\end{subfigure}
\begin{subfigure}[t]{.24\textwidth} 
\centering
\includegraphics[width=1.0\linewidth]{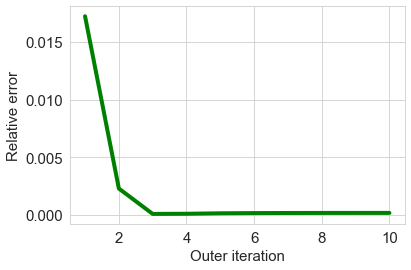}
\caption{\label{fig:lambda3_outer}Relative error in outer iterations}
\end{subfigure}
\caption{Performance of the algorithm when $\lambda_{SE}=3.0$, with true $M=0.75$ and learned $\widehat{M} = 0.736$.}
\end{figure}

	\begin{figure}[H]
\centering
\includegraphics[width=0.5\linewidth]{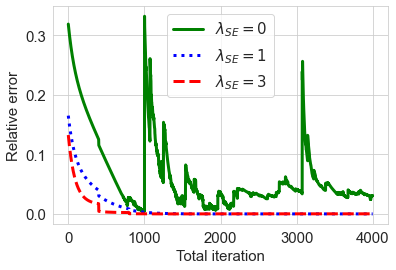}
\caption{\label{fig:experiments_shannon_comparison}Comparison of relative errors with different $\lambda_{SE}$.}
\end{figure}

\paragraph{Results.}
\begin{enumerate}
\item {\it Stability.} As seen from Figure~\ref{fig:experiments_shannon_comparison},
when $\lambda_{SE}=0$, i.e., when there is no exploration,  the algorithm is unstable. Within each outer iteration, the error level fluctuates when the representative agent updates her policy under a fixed mean field information. At the end of each outer iteration, there is a sudden jump in the error when the population updates its mean field policy. In contrast, the algorithm is stable when exploration is included, i.e., when $\lambda_{SE}>0$. 
    \item {\it Speed of  convergence.} As   Figures \ref{fig:lambda1_m} and \ref{fig:lambda3_m} show, Shannon entropy ($\lambda_{SE}>0$)  improves the speed of  convergence to the mean field equilibrium. In fact, the algorithm does not converge without entropy regularization, i.e., when $\lambda_{SE}=0$;  On the other hand, the algorithm converges to the equilibrium solution when $\lambda_{SE}=1$ and $\lambda_{SE}=3$. Moreover, the convergence speed is faster with $\lambda_{SE}=3$ than with $\lambda_{SE}=1$, with the former converging to the mean field equilibrium within three outer iterations and the latter in five outer iterations.
    \item {\it Accuracy of learned mean field equilibrium}. Figures~\ref{fig:lambda1_m} and ~\ref{fig:lambda3_m} show consistency with Theorems \ref{thm:current_shannon} and \ref{thm:current_cross}.
    The algorithm is able to learn the mean field information with small errors ($<5\%$) for both cases $\lambda_{SE}=1$ and $\lambda_{SE}=3$.   
    \item {\it Learning optimal scheduling of the exploration policy.} With given parameters, the variance of the Gaussian mean field policy (a.k.a., the optimal exploration scheduling) is a decreasing function of time $t$ for both  $\lambda_{SE}=1$ and $\lambda_{SE}=3$. Figures \ref{fig:lambda1_V} and \ref{fig:lambda3_V} suggest that  the agent can learn this decreasing function $\{\widehat{\sigma}^2_s\}_{s=0}^T$ with small error $(<5\%)$.
\end{enumerate}

\newpage
\bibliographystyle{plain}
\bibliography{MFG_EE_REF.bib}

\end{document}